\documentclass[letterpaper, 11pt]{article}
\usepackage{amsmath,amsfonts,amsthm,xypic}
\usepackage{hyperref}

\theoremstyle{plain}
\newtheorem{theo}{Theorem}[section]
\newtheorem{cor}[theo]{Corollary}
\newtheorem{lemma}[theo]{Lemma}
\newtheorem{lem}[theo]{Lemma}
\newtheorem{prop}[theo]{Proposition}

\theoremstyle{definition}
\newtheorem{example}[theo]{Example}
\newtheorem{definition}[theo]{Definition}

\newenvironment{renumerate}
{
\begin{enumerate}}
{\end{enumerate}
}

{\vskip6pt%
\noindent%
{\it Remarks}: %
\begin{renumerate}}%
{\end{renumerate}\vskip6pt}
\newenvironment{ex}[1]%
{\begin{example}\label{#1}}%
{\end{example}}

\newcommand{\euler}{\text{\tt e}}

\newcommand{\CC}{\mathbb C}
\newcommand{\R}{\text{${\mathbb R}$}}
\newcommand{\C}{\text{$\mathbb C$}}

\renewcommand{\frak}[1]{\text{$\mathfrak{#1}$}}
\newcommand{\J}{\text{$\mathcal{J}$}}
\newcommand{\JJ}{\text{$\mathcal{J}$}}

\newcommand{\ga}{a}

\newcommand{\e}{\text{$\varepsilon$}}

\newcommand{\w}{\text{$\omega$}}
\newcommand{\gf}{\text{$\varphi$}}

\newcommand{\del}{\text{$\partial$}}

\newcommand{\delbar}{\text{$\overline{\partial}$}}
\newcommand{\tensor}{\otimes}

\newcommand{\mc}[1]{\text{$\mathcal{#1}$}}

\newcommand{\into}{\longrightarrow}
\newcommand{\noqed}{\let\qed\relax}

\newcommand{\Gg}{\mathfrak{g}}

\newcommand{\gO}{\text{$\Omega$}}
\renewcommand{\to}{\longrightarrow}

\newcommand{\gcs}{generalized complex structure}

\newcommand{\gcss}{generalized complex structures}

\newcommand{\gk}{generalized K\"ahler}
\newcommand{\gks}{generalized K\"ahler structure}

\newcommand{\wrt}{with respect to}
\renewcommand{\iff}{if and only if}
\newcommand{\Ann}{\mathrm{Ann}}
\newcommand{\Kperp}{\text{$K^{\perp}$}}
\newcommand{\Diff}{\mathrm{Diff}}

\newcommand{\RR}{\mathbb R}
\newcommand{\LL}{\mathcal L}

\newcommand{\lra}{\longrightarrow}

\newcommand{\dd}{\mathcal{D}}

\newcommand{\BAR}{\overline}
\newcommand{\IP}[1]{\langle #1\rangle}

\newcommand{\eps}{\varepsilon}

\newcommand{\E}{{E}}
\newcommand{\ad}{\mathrm{ad}}

\newcommand{\Kc}{K_{\mathbb{C}}}
\newcommand{\Aa}{\mathfrak{a}}
\newcommand{\Hh}{\mathfrak{h}}

\newcommand{\Cour}[1]{[#1]}

\title{\bf Reduction of Courant algebroids and generalized complex structures}
\author{Henrique Bursztyn\thanks{{\tt
henrique@impa.br}}\ \ \ \  Gil R. Cavalcanti\thanks{{\tt
gilrc@maths.ox.ac.uk}}\ \ \ \ Marco Gualtieri\thanks{{\tt
mgualt@math.mit.edu}} }

\date{}
\begin{document}
\maketitle

\abstract{We present a theory of reduction for Courant algebroids as
well as Dirac structures, generalized complex, and generalized
K\"ahler structures which interpolates between holomorphic reduction
of complex manifolds and symplectic reduction. The enhanced symmetry group of a Courant algebroid leads
us to define \emph{extended} actions and a generalized notion of
moment map.  Key examples of generalized K\"ahler reduced spaces
include new explicit bi-Hermitian metrics on $\CC P^2$. }

\tableofcontents \pagestyle{headings}

\pagebreak
\section{Introduction}
In the presence of a symmetry, a given geometrical structure may,
under suitable conditions, pass to the quotient.  Often, however,
the quotient does not inherit the same type of geometry as the
original space; it may be necessary to pass to a further
\emph{reduction} for this to occur.  For example, a complex manifold
$M$ admitting a holomorphic $S^1$ action certainly does not induce a
complex structure on $M/S^1$; rather, one considers the
complexification of this action to a $\CC^*$ action, whose quotient,
under suitable conditions, inherits a complex structure. Similarly,
the quotient of a symplectic manifold by a symplectic $S^1$ action
is never symplectic; rather it is endowed with a natural Poisson
structure, whose leaves are the symplectic reduced spaces one
desires.

In this paper we consider the reduction of generalized geometrical
structures such as Dirac structures and generalized complex
structures.  These are geometrical structures defined not on the
tangent bundle of a manifold but on the sum $TM\oplus T^*M$ of the
tangent and cotangent bundles (or, more generally, on an exact
Courant algebroid).  These structures interpolate between many of
the classical geometries such as symplectic and Poisson geometry,
the geometry of foliations, and complex geometry.  As a result the
quotient procedure described in this paper interpolates between the
known methods of reduction in these cases.

The main conceptual advance required to understand the reduction
of generalized geometries is the fact that one must extend the
notion of action of a Lie group on a manifold.  Traditional
geometries are defined in terms of the Lie bracket of vector
fields, whose symmetries are given precisely by diffeomorphisms.
As a result, one considers reduction in the presence of a group
homomorphism from a Lie group into the group of diffeomorphisms.
The Courant bracket, on the other hand, has an enhanced symmetry
group which is an abelian extension of a diffeomorphism group by
the group of closed 2-forms. For this reason one must consider
actions which may have components acting nontrivially on the
Courant algebroid while leaving the underlying manifold fixed.  To
formalize this insight, we introduce the notion of a \emph{Courant
algebra}, and explain how it acts on a Courant algebroid in a way
which extends the usual action of a Lie algebra by tangent vector
fields.

A surprising benefit of this point of view is that the concept of
\emph{moment map} in symplectic geometry obtains a new
interpretation as an object which controls the extended part of the
action mentioned above, that is, the part of the action trivially
represented in the diffeomorphism group.

In preparing this article, the authors drew from a wide variety of
sources, all of which provided hints toward the proper framework for
generalized reduction.  First, the literature on holomorphic
reduction of complex manifolds as well as the field of Hamiltonian
reduction of symplectic manifolds in the style of Marsden-Weinstein
\cite{MW}. Also, in the original work of Courant and Weinstein
(\cite{Courant}, \cite{CourWein}) where the Courant bracket is
introduced, some preliminary remarks about quotients can be found;
subsequent formulations of reduction of Dirac structures appear in
\cite{BlaVan,BurCra,LWX98}. Most influential, however, has been the
work of physicists on the problem of finding gauged sigma models
describing supersymmetric sigma models with isometries.  The reason
this is relevant is that the geometry of a general $N=(2,2)$
supersymmetric sigma model is equivalent to generalized K\"ahler
geometry~\cite{Gua03}, and so any insight into how to ``gauge'' or
quotient such a model provides us with guidance for the geometrical
reduction problem.  Our sources for this material have been the work
of Hull, Ro\v{c}ek, de Wit, and
Spence~(\cite{HullRoc},\cite{HullSpence}), Witten~\cite{Witten}, and
Figueroa-O'Farill and Stanciu~\cite{FFS}.  More recently in the
physics literature, the gauging conditions have been re-interpreted
in terms of the Courant bracket~\cite{Fig05}, a point of view which
we develop and expand upon in this paper as well.   Finally, in
recent work of Hitchin~\cite{Hit05}, a natural generalized K\"ahler
structure on the moduli space of instantons on a generalized
K\"ahler 4-manifold is constructed by a method which amounts to an
infinite-dimensional generalized K\"ahler quotient.

The paper is organized as follows.  In Section~\ref{two} we review
the definition of Courant algebroid, describe its group of
symmetries, and define the concept of extended action.  This
involves the definition of a Courant algebra, a particular kind of
Lie 2-algebra. In this section we also define a moment map for an
extended action. In Section~\ref{quotient courant algebroid} we
describe how an extended action on an exact Courant algebroid gives
rise to reduced spaces equipped with induced exact Courant
algebroids. It turns out that, even if the original Courant
algebroid has trivial 3-form curvature, its reduced spaces may have
nontrivial curvature. In Section~\ref{four} we arrive at the
reduction procedure for generalized geometries, introducing an
operation which transports Dirac structures from a Courant algebroid
to its reduced spaces. This operation generalizes both the operation
of Dirac push-forward and pull-back outlined in~\cite{BurRadko}. In
Section~\ref{five} we apply this procedure to reduce generalized
complex structures and provide several examples, including some with
interesting type change. Finally in Section~\ref{six} we study a way
to transport a generalized K\"ahler structure to the reduced spaces.
This is very much in the spirit of the usual K\"ahler reduction
procedure (see~\cite{GuSt82,Kirwan}). Finally we present two examples of
generalized K\"ahler reduction: we produce generalized K\"ahler
structures on $\CC P^2$ with type change, first along a triple line
(an example of which has been found in~\cite{Hit05} using a
different method) and second, along three distinct lines in the
plane.  These examples are particularly significant since they
provide explicit  bi-Hermitian metrics on $\CC P^2$.

Recently there has been a great deal of interest in porting the
techniques of Hamiltonian reduction to the setting of generalized
geometry.  The authors are aware of four other groups who have
worked independently on this topic: Lin and Tolman~\cite{LT05},
Stienon and Xu~\cite{SX05},  Hu~\cite{Hu05}, and
Vaisman~\cite{Va05}.

{\bf Acknowledgements:}  The authors wish to thank M. Crainic, R.
Fernandes, N. Hitchin, C. Hull,  L. Jeffrey, A. Kapustin, Y.
Karshon, F. Kirwan, Y. Li, E. Meinrenken, A. Weinstein, and E.
Witten for many helpful conversations along the way.  We also thank
the Fields Institute, Oxford's Mathematical Institute, IMPA, NSERC
and EPSRC for supporting this project.  Finally, we thank the
referee for many helpful suggestions.

\section{Symmetries of the Courant bracket}\label{two}
In this section we introduce an extended notion of group action on a
manifold preserving twisted Courant brackets. We start by recalling
the definition and basic properties of Courant algebroids.

\subsection{Courant algebroids}\label{courant algebroids}

Courant algebroids were introduced in~\cite{LWX97} in order to
axiomatize the properties of the Courant bracket, an operation on
sections of $TM\oplus T^*M$ extending the Lie bracket of vector
fields.

A {\it Courant algebroid} over a manifold $M$ is a vector bundle $\E
\to M$ equipped with a fibrewise nondegenerate symmetric bilinear
form $\IP{\cdot,\cdot}$, a bilinear bracket $\Cour{\cdot,\cdot}$ on
the smooth sections $\Gamma(E)$, and a bundle map $\pi: \E\to TM$
called the \textit{anchor}, which satisfy the following conditions for all
$e_1,e_2,e_3\in \Gamma(E)$ and $f\in C^{\infty}(M)$:
\begin{itemize}
\item[C1)] $\Cour{e_1,\Cour{e_2,e_3}} = \Cour{\Cour{e_1,e_2},e_3} +
\Cour{e_2,\Cour{e_1,e_3}}$,
\item[C2)] $\pi(\Cour{e_1,e_2})=[\pi(e_1),\pi(e_2)]$,
\item[C3)] $\Cour{e_1,fe_2}=f\Cour{e_1,e_2}+ (\pi(e_1)f) e_2$,
\item[C4)] ${\pi(e_1)}\IP{e_2,e_3}= \IP{\Cour{e_1,e_2},e_3}
+ \IP{e_2, \Cour{e_1, e_3}}$,
\item[C5)] $[e_1,e_1] = \dd\IP{e_1,e_1}$,
\end{itemize}
where $\dd =\tfrac{1}{2}\pi^*\circ d: C^{\infty}(M)\lra \Gamma(E)$
(using $\IP{\cdot,\cdot}$ to identify $\E$ with $\E^*$).

We see from axiom $C5)$ that the bracket is not skew-symmetric, but
rather satisfies
$$
[e_1,e_2]=-[e_2,e_1]+2\mathcal{D}\IP{e_1,e_2}.
$$
Since the left adjoint action is a derivation of the bracket (axiom
$C1)$), the pair $(\Gamma(E),[\cdot,\cdot])$ is a \emph{Leibniz
algebra}~\cite{Loday}. Note that the skew-symmetrization of this
bracket does not satisfy the Jacobi identity; as was shown
in~\cite{RW}, a Courant algebroid provides an example of an
$L_\infty$ algebra. We now briefly describe \v{S}evera's
classification of exact Courant algebroids.

\begin{definition} A Courant algebroid is {\it exact} if
the following sequence is exact:
\begin{equation}\label{exact}
0 \longrightarrow T^*M \stackrel{\pi^*}{\longrightarrow} \E \stackrel{\pi}{\longrightarrow}
TM \longrightarrow 0
\end{equation}
\end{definition}
Given an exact Courant algebroid, we may always choose a right
splitting $\nabla: TM\to \E$ which is \emph{isotropic}, i.e. whose
image in $\E$ is isotropic with respect to $\IP{\cdot,\cdot}$. Such
a splitting has a curvature 3-form $H \in \Omega_{cl}^3(M)$ defined
as follows: for $X,Y\in \Gamma(TM)$,
\begin{equation}\label{eq:curv}
i_Yi_X H = {2}s\Cour{\nabla(X),\nabla (Y)},
\end{equation}
where $s:\E\lra T^*M$ is the induced left splitting.  Using the
bundle isomorphism $\nabla+\tfrac{1}{2}\pi^*:TM\oplus T^*M \to \E$,
we transport the Courant algebroid structure onto $TM\oplus T^*M$.
Given $X+\xi,Y+\eta\in \Gamma(TM\oplus TM^*)$, we obtain for the
bilinear pairing
\begin{equation}\label{eq:pairing}
\IP{X+\xi, Y+\eta} = \frac{1}{2}(\eta(X) + \xi(Y)),
\end{equation}
and the bracket becomes
\begin{equation}\label{eq:Hcour}
\Cour{X+\xi, Y+\eta}_H= [X,Y] + \mc{L}_X \eta - i_Y d\xi  + i_Y i_X
H,
\end{equation}
which is the \emph{$H$-twisted Courant bracket on $TM\oplus T^*M$}
\cite{SW01}.  Isotropic splittings of $\eqref{exact}$ differ by
2-forms $b\in \Omega^2(M)$, and a change of splitting modifies the
curvature $H$ by the exact form $db$. Hence the cohomology class
$[H]\in H^3(M,\mathbb{R})$, called the \textit{\v{S}evera class}, is
independent of the splitting and determines the exact Courant
algebroid structure on $\E$ completely. When this class is integral,
the exact Courant algebroid may be viewed as a generalized Atiyah
sequence associated to a connection on an $S^1$ gerbe. In this
sense, exact Courant algebroids arise naturally from the study of
gerbes.

We now determine the symmetry group of an exact Courant algebroid,
that is, the group of bundle automorphisms preserving the Courant
algebroid structure.
\begin{definition}
The automorphism group $\mathrm{Aut}(E)$ of a Courant algebroid $E$
is the group of bundle automorphisms $F:E\lra E$ covering
diffeomorphisms $\varphi:M\lra M$ such that
\begin{itemize}
\item[i)] $\varphi^*\IP{F\ \cdot,F\ \cdot}=\IP{\cdot,\cdot}$, i.e.
$F$ is orthogonal,
\item[ii)] $[F\ \cdot,F\ \cdot]=F[\cdot,\cdot]$, i.e.
$F$ is bracket-preserving,
\item[iii)] $\pi\circ F = \varphi_*\circ\pi$, i.e. $F$ is compatible with the anchor.
\end{itemize}
One can easily verify, using axiom $\mathrm{C}3)$, that
compatibility with the anchor is implied by requirements
$\mathrm{i)}$ and $\mathrm{ii)}$ (See~\cite{Gua03,Wald}).

Similarly, the Lie algebra of derivations $\mathrm{Der}(E)$ is the Lie
algebra of linear first order differential operators $D_X$ on
$\Gamma(E)$, covering vector fields $X\in\Gamma(TM)$ such that
$X\IP{\cdot,\cdot} = \IP{D_X\ \cdot,\cdot}+\IP{\cdot,D_X\ \cdot}$
and $D_X[\cdot,\cdot]=[D_X\ \cdot,\cdot]+[\cdot,D_X\ \cdot]$.
\end{definition}

In the case of an exact Courant algebroid, one may choose an
isotropic splitting of the anchor, inducing an isomorphism $E\cong
TM\oplus T^*M$ as above, with bilinear pairing given
by~\eqref{eq:pairing} and bracket given by~\eqref{eq:Hcour}. Now
suppose that $F\in \mathrm{Aut}(E)$ covers $\varphi\in\mathrm{Diff}(M)$.
Note that $\varphi$ lifts naturally to $\Phi=\varphi_*+
(\varphi^*)^{-1}\in \mathrm{End}(TM\oplus T^*M)$, which satisfies
$$[\Phi\ \cdot,\Phi\ \cdot]_H = \Phi[\cdot,\cdot]_{\varphi^*H}.$$
Therefore $\Phi^{-1}F$ is a fiber-preserving orthogonal map on
$TM\oplus T^*M$ compatible with the anchor, which implies that
it must be the orthogonal action of a 2-form $B\in\Omega^2(M)$ via
$e^B:X+\xi\mapsto X+\xi+i_XB$ \cite{Gua03}.  Since these ``gauge
transformations'' satisfy
\[
[e^B\ \cdot, e^B\ \cdot]_H = e^B[\cdot,\cdot]_{H+dB},
\]
we see that $F=\Phi e^B$ is an automorphism if and only if
$H-\varphi^*H = dB$.  Therefore, the automorphism group consists of
ordered pairs $(\varphi,B)\in \Diff(M)\times \Omega^2(M)$ such that
$H-\varphi^*H = dB$, giving rise to the following
splitting-independent description.
\begin{prop}
The automorphism group of an exact Courant algebroid $E$ is an
extension of the diffeomorphisms preserving the cohomology class
$[H]$ by the abelian group of closed 2-forms:
\begin{equation}\label{extenex}
\xymatrix{0\ar[r]&\Omega^2_{cl}(M)\ar[r]&\mathrm{Aut}(E)\ar[r]&\Diff_{[H]}(M)\ar[r]&0}
\end{equation}
If $M$ is compact, the extension class in group cohomology is
represented by the cocycle
\[
c(\varphi_1,\varphi_2) =
{\varphi_1^*}^{-1}(Q-{\varphi_2^*}^{-1}Q\varphi_2^*)(H-\varphi_1^*H),
\]
where $Q=d^*G$, and $d^*,G$ are the codifferential and Green
operator with respect to a Riemannian metric.
\end{prop}
\begin{proof}
Given an isotropic splitting of $E$ with curvature $H$, and a
Riemannian metric on the compact manifold $M$, we split the
sequence~\eqref{extenex} via the map $s:\Diff_{[H]}(M)\lra
\mathrm{Aut}(E)$ given by $s(\varphi) = (\varphi, B_\varphi)$, where
$B_\varphi = Q(H-\varphi^*H)$.  One can easily verify that
$dB_\varphi = H-\varphi^*H$ and that
$$s(\varphi_1)s(\varphi_2)(s(\varphi_1\varphi_2))^{-1} =
(1,c(\varphi_1,\varphi_2)),$$ yielding the extension class.
\end{proof}

Differentiating a 1-parameter family of automorphisms $F_t =
\Phi_te^{tB}$, $F_0=\mathrm{Id}$, we see that the Lie algebra
$\mathrm{Der}(E)$ for a split exact Courant algebroid consists of
pairs $(X,B)\in \Gamma(TM)\oplus \Omega^2(M)$ such that $\LL_X H =
dB$, which act via
\begin{equation}\label{ac}
(X,B)\cdot(Y+\eta) = \LL_X(Y+\eta) + i_YB.
\end{equation}
We then have the following invariant description of derivations.
\begin{prop}
The Lie algebra of infinitesimal symmetries of an exact Courant
algebroid $E$ is an abelian extension of the Lie algebra of smooth
vector fields by the closed 2-forms:
\begin{equation}\label{extenlie}
\xymatrix{0\ar[r]&\Omega^2_{cl}(M)\ar[r]&\mathrm{Der}(E)\ar[r]&\Gamma(TM)\ar[r]&0}
\end{equation}
The extension class in Lie algebra cohomology is represented by the
cocycle
\[
c(X,Y) = d i_X i_Y H.
\]
\end{prop}
\begin{proof}
Given an isotropic splitting of $E$ with curvature $H$, we split
sequence~\eqref{extenlie} via the map $s:\Gamma(TM)\lra
\mathrm{Der}(E)$ given by $s(X)=(X,i_XH)$. Then, using~\eqref{ac},
we find that
\[
[s(X),s(Y)]-s([X,Y]) = (0,c(X,Y)),
\]
as required.
\end{proof}
\noindent This Lie algebra cocycle was also obtained by Hu
in~\cite{Hu05}, where more details can be found.

It is immediately clear from axioms $\mathrm{C1), C4)}$ that
$\Gamma(E)$ acts on itself by derivations via the left adjoint
action $\ad_v(w):=[v,w]$.  Unlike, however, the usual adjoint action
of vector fields on the tangent bundle, the map $\ad: \Gamma(E)\to
\mathrm{Der}(E)$ is neither surjective nor injective, as we now
verify for exact Courant algebroids.
\begin{prop}
Let $E$ be an exact Courant algebroid. Then the adjoint action
$\ad:v\mapsto[v,\cdot]$ induces the following exact sequence:
\begin{equation*}
\xymatrix{0\ar[r]&\Omega^1_{cl}(M)\ar[r]^{\pi^*}&\Gamma(\E)\ar[r]^{\ad}&\mathrm{Der}(E)\ar[r]^\chi&H^2(M,\RR)\ar[r]&0}
\end{equation*}
\end{prop}
\begin{proof}
Given an isotropic splitting, we see from~\eqref{eq:Hcour} that the
kernel of the adjoint action is the space of closed 1-forms.  Given
any derivation $D_X=(X,B)$, we define $\chi(D_X)=[i_XH-B]\in
H^2(M,\R)$, which is surjective by the freedom to choose $B$, and
whose kernel consists of $(X,B)$ such that $B=i_XH-d\xi$, i.e. such
that $(X,B)\cdot (Y+\eta)= [X+\xi,Y+\eta]$, proving exactness.
\end{proof}

\subsection{Extended actions}

Let a Lie group $G$ act on a manifold $M$, so that we have the Lie
algebra homomorphism $\psi:\Gg\lra \Gamma(TM)$. We wish to extend
this action to a Courant algebroid $\E$, making $E$ into a
$G$-equivariant vector bundle, in such a way that the Courant
algebroid structure is preserved.  In this section we show how this
can be done by choosing an extension of $\Gg$ equipped with a
\emph{Courant algebra} structure, and choosing a homomorphism from
this extension to the Courant algebroid $\E$.

\begin{definition}
A \emph{Courant algebra} over the Lie algebra $\Gg$ is a vector
space $\Aa$ equipped with a bilinear bracket
$[\cdot,\cdot]:\Aa\times\Aa\lra \Aa$ and a map $\pi:\Aa\lra\Gg$,
which satisfy the following conditions for all $a_1, a_2,
a_3\in\Aa$:
\begin{itemize}
\item[c1)] $[a_1,[a_2,a_3]]=[[a_1,a_2],a_3]+[a_2,[a_1,a_3]]$,
\item[c2)] $\pi([a_1,a_2])=[\pi(a_1),\pi(a_2)]$.
\end{itemize}
In other words, $\Aa$ is a Leibniz algebra with a homomorphism to
$\Gg$.
\end{definition}
A Courant algebroid provides an example of a Courant algebra over
$\Gg=\Gamma(TM)$, taking $\Aa=\Gamma(E)$.  Using the argument of
Roytenberg-Weinstein~\cite{RW}, one sees that any Courant algebra is
actually an example of a 2-term $L_\infty$-algebra \cite{CrBa,Sta}.

\begin{definition}
An \emph{exact} Courant algebra is one for which $\pi$ is surjective
and $\Hh=\ker\pi$ is abelian, i.e. $[h_1,h_2]=0$ for all
$h_1,h_2\in\Hh$.
\end{definition}
For an exact Courant algebra, one obtains immediately an action of
$\Gg$ on $\Hh$: $g\in\Gg$ acts on $h\in\Hh$ via $g\cdot h = [a, h]$,
for any $a$ such that $\pi(a)=g$. This is well defined since $\Hh$
is abelian, and determines an action by axiom $\mathrm{c1)}$.  In
fact there is a natural exact Courant algebra associated with any
$\Gg$-module, as we now explain.

\begin{example}[Hemisemidirect product]\label{Courrep}
Let $\Gg$ be a Lie algebra acting on the vector space $\Hh$.  Then
$\Aa=\Gg\oplus\Hh$ becomes a Courant algebra over $\Gg$ via the
bracket
\begin{equation}\label{dsb}
[(g_1,h_1),(g_2,h_2)]=([g_1,g_2],g_1\cdot h_2),
\end{equation}
where $g\cdot h$ denotes the $\Gg$-action.  This bracket appeared
in~\cite{KinWein}, where it was called the
\emph{hemisemidirect} product of $\Gg$ with $\Hh$. Note that
in~\cite{WeinOmni}, Weinstein studied the case where
$\Gg=\mathfrak{gl}(V)$ and $\Hh=V$, and called it an \emph{omni-Lie
algebra} due to the fact that, when $\dim V=n$,  any $n$-dimensional
Lie algebra can be embedded in it as a subalgebra.
\end{example}

\begin{definition}[Extended action]
Let $G$ be a connected Lie group acting on a manifold $M$ with
infinitesimal action $\psi:\Gg\lra \Gamma(TM)$.  An \textit{extension} of
this action to a Courant algebroid $\E$ over $M$ is an exact Courant
algebra $\Aa$ over $\Gg$ together with a Courant algebra morphism
$\rho:\Aa\lra \Gamma(\E)$:
\begin{equation*}
\xymatrix{0\ar[r]&\Hh\ar[r]&\Aa\ar[r]\ar[d]_{\rho}&\Gg\ar[d]^{\psi}\ar[r]&0\\
 & & \Gamma(\E)\ar[r]&\Gamma(TM)&}
\end{equation*}
which is such that $\Hh$ acts trivially, i.e.
$(\ad\circ\rho)(\Hh)=0$, and the induced action of $\Gg=\Aa/\Hh$ on
$\Gamma(E)$ integrates to a $G$-action on the total space of $\E$.
\end{definition}

The space of $G$-equivariant functions $f:M\lra \Gg^*$ acts on the
space of extensions of an action $\psi$ by the Courant algebra
$\Aa$, via $\rho(\cdot)\mapsto \rho(\cdot) + \dd \IP{f,\pi(\cdot)}$, where $\dd = \frac{1}{2} \pi^* \circ d$, as above;
this generates an equivalence relation among extensions of actions.
\begin{definition}\label{eqaction}
Extensions $\rho,\rho'$ of a given $G$-action to a Courant algebroid
$E$, for a fixed Courant algebra $\Aa$, are said to be
\emph{equivalent} if they agree upon restriction to $\Hh$ and differ
by a $G$-equivariant function $f:M\lra \Gg^*$, i.e.
\[
\rho'(a)-\rho(a) = \dd\IP{f,\pi(a)}.
\]
\end{definition}

Suppose now that the Courant algebroid in question is exact, as it
will be in many cases of interest. Then an extended action is a
commutative diagram
\begin{equation*}
\xymatrix{0\ar[r]&\Hh\ar[r]\ar[d]^{\nu}&\Aa\ar[r]\ar[d]^{\rho}&\Gg\ar[d]^{\psi}\ar[r]&0\\
0\ar[r]& \Gamma(T^*M)\ar[r]&
\Gamma(E)\ar[r]&\Gamma(TM)\ar[r]&0}
\end{equation*}
such that $\Hh$ acts trivially, which occurs precisely when it
acts via {\it closed} 1-forms, i.e. $\nu(\Hh)\subset
\Omega^1_{cl}(M)$. Furthermore the induced $\Gg$-action on $\E$
must integrate to a $G$-action (a priori, one has only the action
of the universal cover of $G$). In order to make this condition
more concrete, we observe that since we already know that the
$\Gg$-action on $TM$ integrates to a $G$-action, one needs only to
find a $\Gg$-invariant splitting of $\E$ to guarantee that it is a
$G$-bundle, as the splitting $\E=TM\oplus T^*M$ carries a canonical
$G$-equivariant structure.
\begin{prop}\label{avgspl}
Let the  Lie group $G$ act on the manifold $M$, and let
$\xymatrix{\Aa\ar[r]^\pi&\Gg}$ be an exact Courant algebra with a
morphism $\rho$ to an exact Courant algebroid $\E$ over $M$ such
that $\nu(\Hh)\subset\Omega^1_{cl}(M)$.

If $\E$ has a $\Gg$-invariant splitting, then the $\Gg$-action on
$E$ integrates to an action of $G$, and hence $\rho$ is an extended
action of $G$ on $E$.  Conversely, if $G$ is compact and $\rho$ is
an extended action, then by averaging splittings one can always find
a $\Gg$-invariant splitting of $E$.
\end{prop}

The condition that a splitting is $\Gg$-invariant can be expressed
more concretely as follows.  As shown in Section~\ref{courant
algebroids}, a split exact Courant algebroid is isomorphic to the
direct sum $TM\oplus T^*M$, equipped with the $H$-twisted Courant
bracket for a closed 3-form $H$.  In this splitting, therefore, for
each $\ga\in\Aa$ the section $\rho(a)$ decomposes as $\rho(\ga) =
X_\ga + \xi_\ga$, and it acts via $[X_\ga+\xi_\ga,Y + \eta] =
[X_\ga,Y] + \mc{L}_{X_\ga}\eta - i_Y d\xi_\ga + i_Yi_{X_\ga}H$, or
as a matrix,
\begin{equation*}
ad_{\rho(\ga)} =\begin{pmatrix}
\mc{L}_{X_\ga} & 0\\
i_{X_{\ga}}H - d\xi_{\ga} & \mc{L}_{X_{\ga}}
\end{pmatrix}
\end{equation*}
We see immediately from this that the splitting is preserved by this
action if and only if for each $a\in\Aa$,
\begin{equation}\label{equi}
i_{X_{\ga}}H - d\xi_{\ga} =0.
\end{equation}

We now provide a complete description, assuming $G$ to be compact
and $E$ exact, of the simplest kind of extended action, namely one
for which $\Aa = \Gg$.
\begin{definition}\label{def:trivially}
A \emph{trivially} extended $G$-action is one for which $\Aa=\Gg$
and $\pi:\Aa\lra\Gg$ is the identity map, as described by the
commutative diagram:
\begin{equation*}
\xymatrix{\Gg\ar[r]^{\text{id}}\ar[d]^{\rho}&\Gg\ar[d]^{\psi}\\
\Gamma(E)\ar[r]&\Gamma(TM)}
\end{equation*}
\end{definition}

Suppose that $G$ is compact and $E$ exact.  By
Proposition~\ref{avgspl}, we can always find a $\Gg$-invariant
splitting of $E$, so finding a trivially extended action $\rho$ is
equivalent to finding 1-forms $\xi_\ga$ such that $\rho:\ga\mapsto
X_{\ga}+\xi_\ga$ is a Courant algebra homomorphism (here $X_\ga =
\psi(\ga)$ and $\ga\in\Gg$). Preserving the bracket yields
\begin{equation}\label{actex}
\xi_{[a,b]}=\LL_{X_a}\xi_b -i_{X_b}d\xi_a +
i_{X_b}i_{X_a}H=\LL_{X_a}\xi_b,
\end{equation}
where~\eqref{equi} was used in the final equality.
Equation~\eqref{actex} states that $\xi_a$ is an equivariant form,
and condition~\eqref{equi} can be phrased in terms of the Cartan
model for $G$-equivariant cohomology. Recall that the Cartan complex
of equivariant forms is the algebra of equivariant polynomial
functions $\Phi:\Gg\lra \Omega^\bullet(M)$:
\begin{equation*}
\Omega_G^k(M)=\bigoplus_{2p+q=k}(S^p\Gg^*\otimes\Omega^q(M))^G,
\end{equation*}
and the equivariant derivative $d_G$ is defined by
\begin{equation*}
(d_G\Phi)(a)= d(\Phi(a))-i_{X_a}\Phi(a)\ \ \ \forall a\in\Gg.
\end{equation*}
Now consider the form $\Phi(a)= H+\xi_\ga$.  Since the splitting is
$G$-invariant, we have $\LL_{X_\ga}H=0$.  Therefore $\Phi$ is an
equivariant 3-form in the Cartan complex. Computing $d_G\Phi$, we
obtain
\begin{equation*}
d_G\Phi(a)= -\IP{X_a+\xi_a,X_a+\xi_a} = -\IP{\rho(a),\rho(a)}.
\end{equation*}
This shows that the quadratic form $c(a)=\IP{\rho(a),\rho(a)}$ is,
firstly, constant along $M$, but also that it defines an invariant
quadratic form on the Lie algebra $\Gg$, and furthermore,
one which is exact in the $G$-equivariant cohomology of $M$.  If
$c \equiv 0$, i.e. if the action is \emph{isotropic}, then we see that
the existence of a trivially extended $G$-action on $E$ is
determined by the equivariant extension of $[H]$.  This last
condition is well-known to physicists in the context of gauging
sigma models with Wess-Zumino term~\cite{HullSpence}.

\begin{theo}\label{eqphic}
Let $G$ be a compact Lie group. Then trivially extended $G$-actions
on a fixed exact Courant algebroid with prescribed quadratic form
$c(a)=\IP{\rho(a),\rho(a)}$ are, up to equivalence, in bijection
with solutions to $d_G\Phi=c$ modulo $d_G$-exact forms, where
$\Phi(a)=H+\xi_{a}$ is an equivariant 3-form and $[H]\in H^3(M,\RR)$
is the \v{S}evera class of the Courant algebroid.
\end{theo}
\begin{proof}
The equivariant 3-form $\Phi(a)=H+\xi_a$ representing a trivially
extended $G$-action $\rho$ depends on a choice of $\Gg$-invariant
splitting for $E$.  Changing the splitting by a gauge transformation
$b\in\Omega^2(M),\ \LL_{X_a}b=0$, the 3-form changes to
$\Phi(a)=H+db + \xi_a + i_{X_a}b$.  Also, an equivalent extended
action $\rho'$ satisfies $\rho'(a)-\rho(a)=df_a$ for a
$G$-equivariant function $f$.  The resulting equivariant 3-form is
\[
\Phi'(a) = H + \xi_a + (db+i_{X_a}b + df_a).
\]
But this is precisely the addition to $\Phi$ of an equivariantly
exact 3-form, i.e. $\Phi' = \Phi + d_G\beta$, $\beta(a) = b + f_a$,
proving the result.
\end{proof}

\subsection{Moment maps for extended actions}
Suppose that we have an extended $G$-action on an exact Courant
algebroid as in the previous section,  so that we have the map
$\nu:\Hh\lra \Omega^1_{cl}(M)$. Because the action is a Courant
algebra morphism, this map is $\Gg$-equivariant in the sense
\begin{equation}\label{eqmom}
\nu(g\cdot h) = \LL_{\psi(g)}\nu(h).
\end{equation}
Therefore we are led naturally to the definition of a moment map for
this extended action, as an equivariant factorization of $\nu$
through the smooth functions.

\begin{definition}
A {\it moment map} for an extended $\Gg$-action on an exact
Courant algebroid is a $\Gg$-equivariant map $\mu:\Hh\lra
C^{\infty}(M,\RR)$ satisfying $d \mu = \nu$, i.e. such that
the following diagram commutes:
\begin{equation*}
\xymatrix{ &
\Hh\ar[dl]_\mu\ar[d]^\nu\\C^{\infty}(M)\ar[r]^{d}&\Gamma(T^*M)}
\end{equation*}
Note that $\mu$ may be alternatively viewed as an equivariant map
$\mu:M\lra \Hh^*$.
\end{definition}
A moment map can be found only if two obstructions vanish. The
first one is the induced map to cohomology $\nu_*:\Hh\lra
H^1(M,\RR)$. Since \eqref{eqmom} implies that $\nu_*$ always
vanishes on $\Gg\cdot\Hh\subset\Hh$, the first obstruction may be
defined as an element
\begin{equation*}
o_1\in H^0(\Gg,\Hh^*)\otimes H^1(M,\RR),
\end{equation*}
where the first term denotes Lie algebra cohomology with values in
the module $\Hh^*$.  When this obstruction vanishes we may choose a
lift $\tilde\mu:\Hh\lra C^{\infty}(M)$.  The second obstruction
results from the failure of this lift to be equivariant: consider
the quantity $c(g,h)=\tilde\mu(g\cdot h) -\LL_{\psi(g)}\tilde\mu(h)$
for $g\in\Gg,\ h\in\Hh$.  From~(\ref{eqmom}) we conclude that $c$ is
a constant function along $M$.  It is easily shown that this
discrepancy, modulo changes of lift, defines an obstruction class
\begin{equation*}
o_2\in H^1(\Gg,\Hh^*).
\end{equation*}
\begin{prop}
A moment map for an extended $\Gg$-action exists if and only if the
obstructions $o_1\in H^0(\Gg,\Hh^*)\otimes H^1(M,\RR)$ and $o_2\in
H^1(\Gg,\Hh^*)$ vanish.  When it exists, a moment map is unique up
to the addition of an element $\lambda\in \Ann(\Gg\cdot\Hh)\subset
\Hh^*$.
\end{prop}

We now show how the usual notions of symplectic and Hamiltonian
actions fit into the framework of extended actions of Courant
algebras.
\begin{example}[Symplectic actions]\label{sym}
Let $G$ be a Lie group acting on a symplectic manifold $(M,\omega)$
preserving the symplectic form, and let $\psi:\Gg\lra \Gamma(TM)$
denote the infinitesimal action.  We now show that there is a
natural extended action of the hemisemidirect product Courant
algebra $\Gg\oplus\Gg$ on the standard Courant algebroid $TM\oplus
T^*M$ with $H=0$. As described in Example~\ref{Courrep}, the Courant
algebra is described by the sequence
\begin{equation*}
\xymatrix{0\ar[r]&\Gg\ar[r]&\Gg\oplus\Gg\ar[r]^{\ \
\pi}&\Gg\ar[r]&0}
\end{equation*}
and is equipped with the bracket
\begin{equation}\label{adj1}
[(g_1,h_1),(g_2,h_2)]=([g_1,g_2],[g_1,h_2]),
\end{equation}
Then define the action $\rho:\Gg\oplus\Gg\lra \Gamma(TM\oplus T^*M)$
by
\begin{equation*}
\rho(g,h) = X_g + i_{X_h}\omega,
\end{equation*}
where $X_g=\psi(g)$, for $g\in\Gg$, and $\omega$ is the symplectic
form.  Then since
\begin{equation*}
[X_{g_1}+i_{X_{h_1}}\omega,X_{g_2}+i_{X_{h_2}}\omega] =
[X_{g_1},X_{g_2}] + \LL_{X_{g_1}}i_{X_{h_2}}\omega = X_{[g_1,g_2]} +
i_{X_{[g_1,h_2]}}\omega,
\end{equation*}
we see that $\rho$ is a Courant morphism.

The question of finding a moment map for this extended action then
becomes one of finding an equivariant map $\mu:\Gg\lra C^{\infty}(M)$
such that
\begin{equation*}
d(\mu_g) = i_{X_g}\omega.
\end{equation*}
Hence we recover the usual moment map for a Hamiltonian action on a
symplectic manifold.
\end{example}

Note that in this formalism, the notion of moment map is no longer
tied to the geometry, i.e. the symplectic form. Instead, it is a
constituent of the extended action. In fact, given an equivariant
map $\mu:M\lra \Hh^*$ for a $\Gg$-module $\Hh$, one can naturally
construct an extended action for which $\mu$ is a moment map, as we
now indicate.
\begin{prop}\label{prop:will need later}
Given a $\Gg$-equivariant map $\mu:M\lra\Hh^*$, where $M$ is a
$G$-space and $\Hh$ a $\Gg$-module, there is an induced extended
action of the Courant algebra $\Gg\oplus\Hh$ with
bracket~(\ref{dsb}) on the exact Courant algebroid $TM\oplus T^*M$
with $H=0$, given by
\begin{equation*}
\rho:(g,h)\mapsto X_g + d(\mu_h),
\end{equation*}
where as before $X_g=\psi(g)$ is the infinitesimal $\Gg$-action.
\end{prop}
More generally, given a trivially extended action $\rho:\Gg\into
\Gamma(E)$ on an exact Courant algebroid, it can be extended to an
action of $\Gg\oplus\Hh$ as above by any equivariant map $\mu:M\lra
\Hh^*$ via the same formula
\begin{equation*}
\tilde\rho:(g,h)\mapsto \rho(g) + d(\mu_h).
\end{equation*}


\section{Reduction of Courant algebroids}\label{quotient courant algebroid}

In this section we develop a reduction procedure for exact Courant
algebroids which can be seen as an ``odd'' analog of the usual
notion of symplectic reduction due to Marsden and Weinstein
\cite{MW}. A key observation is that an extended $G$-action on an
exact Courant algebroid $E$ over a manifold $M$ does \emph{not}
necessarily induce an exact Courant algebroid on $M/G$, but rather
one may need to pass to a suitably chosen submanifold $P\subset
M$, in such a way that the \emph{reduced space} $P/G=M_{red}$
obtains an exact Courant algebroid. This is directly analogous to
the well-known fact that, for a symplectic $G$-space $M$, the
reduced spaces are the leaves of the Poisson structure inherited
by $M/G$.

\subsection{Reduction procedure}\label{subsec:redCA}

As we saw in the previous section, an extended action of a
connected Lie group $G$ on a Courant algebroid $E$ over $M$ makes
$E$ into an equivariant $G$-bundle in such a way that the Courant
structure is preserved by the $G$-action. Therefore, assuming the
$G$-action on the base to be free and proper, we obtain a Courant
algebroid $E/G$ over $M/G$. However, even if $E$ were exact, $E/G$
would certainly \textit{not} be an exact Courant algebroid, since
its rank is too large. We will see in this section how this
construction can be modified so as to yield an \textit{exact}
reduced Courant algebroid.

So let $G$ be a connected Lie group and $E$ be an \textit{exact}
Courant algebroid over $M$. The first basic observation is that an
extended action $\rho:\Aa\times M\lra E$ determines two natural
distributions in $E$: the image of $\rho$, $K=\rho(\Aa)$, and its
orthogonal, $K^\perp$. Recall that the action of $g\in\Gg$ on any
generating section $\rho(a)$ of $K$ is simply
\begin{equation*}
g \cdot \rho(a) = [\rho(\tilde g), \rho(a)] = \rho([\tilde g, a])
\in \Gamma(K),
\end{equation*}
where $\tilde{g}\in \Aa$ is any lift  of $g$, $\pi(\tilde g)=g$.
It follows that $K$ is a $G$-invariant distribution and, since the
$G$-action on $E$ preserves the symmetric pairing
$\IP{\cdot,\cdot}$, $K^\perp$ is a $G$-invariant distribution as
well. It is clear that if $\rho$ has constant rank, then the
distributions $K$ and $K^\perp$ are subbundles of $E$, but we do
not make this global assumption at this point.

\begin{definition}
Given an extended action with image distribution
$\rho(\Aa)=K\subset E$, define the \emph{big distribution}
$\Delta_b = \pi(K+K^\perp)\subset TM$ and the \emph{small
distribution} $\Delta_s = \pi(K^\perp)\subset TM$.  These are
$G$-invariant distributions.
\end{definition}
For the construction of \textit{reduced spaces}, we will need to
consider submanifolds tangent to these distributions. The integrability problem for
$\Delta_s$ and $\Delta_b$ will be discussed in Proposition \ref{lem:consrank},
but we make a few observations now. First, note that  $\Delta_s$
satisfies
\begin{equation}\label{kerhh}
\Delta_s  = \Ann (\rho(\Hh)).
\end{equation}
Since the space of sections of $\rho(\Hh)$ is generated by closed 1-forms, it follows that
$\Delta_s$ is an integrable distribution around the points where $\rho(\Hh)$ has locally constant rank.
As we shall see in Section \ref{subsec:severaclass}, in
the presence of a moment map $\mu:M\lra \Hh^*$, $\Delta_s$
coincides with the distribution tangent to the level sets, whereas
$\Delta_b$ is the distribution tangent to the $G$-orbits of the
level sets.

In general, since $\pi(K)$ is the distribution tangent to the
$G$-orbits on $M$, the $G$-orbit of any leaf of $\Delta_s$ (if
smooth) is then a leaf of $\Delta_b$. (Here a \emph{leaf} of a
distribution is taken to mean a maximal connected integral
submanifold.) In particular, any leaf of $\Delta_b$ is
$G$-invariant. These observations allow us to prove the following
useful lemma.

\begin{lemma}\label{cstrk}
Let $P\subset M$ be a leaf of the big distribution $\Delta_b$ on
which $G$ acts freely and properly, and suppose $\rho(\Hh)$ has
constant rank along $P$. Then $K$ and $K\cap K^\perp$ both have
constant rank along $P$.
\end{lemma}

\begin{proof}
Since $G$ acts freely on $P$, $\pi(K)=\psi(\Gg)$ has constant rank
along $P$.  Further, as $\rho(\Hh)$ also has constant rank along $P$, it
follows that $\rho(\Aa)=K$ has constant rank along $P$.

From \eqref{kerhh} and the discussion following it, we conclude that $\Delta_s|_P\subset TP$
is integrable, defining a regular foliation in $P$. Moreover, $P$
is the $G$-orbit of a leaf $S$ of $\Delta_s|_P$.

On the other hand, because $\rho$ is a Courant morphism, we have
for all $a\in\Aa$,
\begin{equation}\label{cab}
\rho([a,a]) = [\rho(a),\rho(a)]=\mathcal{D}\IP{\rho(a),\rho(a)}.
\end{equation}
Since $[a,a]\in \Hh$, it follows that $\rho([a,a])|_{TS}=0$, so we
see that $\IP{\rho(a),\rho(a)}$ is constant along $S$. Hence we obtain an
induced inner product on $\Aa$, $(a,b)\mapsto
\IP{\rho(a),\rho(b)}|_S$, whose null space, modulo
$\ker\rho|_\Hh$, maps isomorphically onto $K\cap K^\perp$. Hence
$K\cap K^\perp$ has constant rank along $S$. But $K\cap K^\perp$
is $G$-invariant, so it must have constant rank over the entire
big leaf $P$.
\end{proof}

Hence, under the assumptions of Lemma \ref{cstrk}, $K$ and $K\cap
K^\perp$ are $G$-invariant vector bundles over $P$, so we can
consider the quotient vector bundle
\begin{equation}\label{eq:Ereduc}
E_{red} = \left. \frac{K^\perp|_P}{K\cap K^\perp|_P}\right/G
\end{equation}
over $M_{red}:=P/G$, the \textit{reduced space}. It is clear that
$E_{red}$ inherits a nondegenerate symmetric pairing from the one
in $E$. The next theorem shows that $E_{red}$ carries in fact a
Courant algebroid structure. We call it the \textit{reduced
Courant algebroid}.

\begin{theo}\label{courred}
Let $E$ be an exact Courant algebroid over $M$ and $\rho:\Aa\to
\Gamma(E)$ be an extended $G$-action.
 Let $P\subset M$ be a leaf of $\Delta_b$ on which $G$ acts
freely and properly, and over which $\rho(\Hh)$ has constant rank.
Then the Courant bracket on $E$ descends to $E_{red}$ and makes it
into a Courant algebroid over $M_{red}=P/G$ with surjective
anchor. If $K$ is isotropic then $E_{red}$ is an exact Courant
algebroid; in general, it is exact if and only if the following
holds along $P$:
\begin{equation}\label{ex}
\pi(K)\cap\pi(K^\perp) = \pi(K\cap K^\perp).
\end{equation}
\end{theo}

\begin{proof}
By Lemma~\ref{cstrk}, $K^\perp$ and $K\cap K^\perp$ are
$G$-invariant bundles over $P$, and hence $E_{red}$ is a vector
bundle over $M_{red}=P/G$ equipped with a nondegenerate symmetric
pairing. We will now check that $E_{red}$ inherits a Courant
bracket.

Let $\tilde v, \tilde w\in \Gamma(E)$ be extensions of
$G$-invariant sections of $K^\perp$ over $P$. Note that the
bracket $\Cour{\tilde v,\tilde w}$ restricted to $P$ is a section
of $K^\perp|_P$: for $a \in \Aa$, by C4) we have
\begin{eqnarray*}
\IP{\rho(a),\Cour{\tilde v,\tilde w}}&=&
-\IP{[\tilde{v},\rho(a)],\tilde{w}} +
{\pi(\tilde{v})}\IP{\rho(a),\tilde{w}}\\
&=& \IP{[\rho(a),\tilde{v}],\tilde{w}}
+{\pi(\tilde{w})}\IP{\rho(a),\tilde{v}} +
{\pi(\tilde{v})}\IP{\rho(a),\tilde{w}},
\end{eqnarray*}
which vanishes along $P$ since $\tilde v$ is an invariant section
of $K^\perp$ there, and $\IP{\rho(a),\tilde{v}}|_P \equiv
\IP{\rho(a),\tilde{w}}|_P\equiv 0$. As a result, $\Cour{\tilde
v,\tilde w}|_P$ is again a $G$-invariant section of $K^\perp|_P$.

To describe the dependence of $[\tilde{v},\tilde{w}]|_P$ with
respect to the extensions chosen, consider a section of $E$
vanishing along $P$, i.e., a section of the form $sf$, where $s
\in \Gamma(E)$ and $f\in C^\infty(M,\mathbb{R})$, with $f|_P\equiv
0$. Then $ \Cour{\tilde{v},sf}=f\Cour{\tilde v, s} + (\pi(\tilde
v)f)s,$ which vanishes upon restriction to $P$, since $f|_P\equiv
0$ and $\pi(\tilde v)$ is tangent to $P$ there. Therefore $
\Cour{\tilde v,sf}|_P\equiv 0$. On the other hand, since
$\Cour{sf,\tilde{w}}= -\Cour{\tilde{w},fs}+ \pi^*d(\IP{s,\tilde
w}f),$ it follows that
\begin{equation}\label{eq:df}
\Cour{sf,\tilde w}|_P=\IP{s,\tilde w}\pi^*df|_P.
\end{equation}
But $df|_P\in\Ann(TP)$ and $\Ann(TP)=\Ann(\pi(K+K^\perp)|_P)
=\{\xi\in T^*M|_P\;|\; \pi^*\xi \in K\cap K^\perp|_P \},$ so the
right-hand side of \eqref{eq:df} is a section of $K\cap
K^\perp|_P$. It follows that if $\tilde v, \tilde w, \tilde v',
\tilde w'$ are sections of $E$ extending $G$-invariant sections of
$K^\perp|_P$ such that $(\tilde v'-\tilde v)|_P = 0$ and $(\tilde
w'-\tilde w)|_P = 0$, then
$$
(\Cour{\tilde v,\tilde w}-\Cour{\tilde v',\tilde w'})|_P\in
\Gamma(K\cap K^\perp|_P)^G.
$$
Hence the bracket on invariant sections of $K^\perp|_P$ is
well-defined modulo sections of $K\cap K^\perp|_P$, which means
that we have a well-defined bracket
$$
\Gamma(K^\perp|_P)^G \times \Gamma(K^\perp|_P)^G\to
\Gamma(E_{red}) = \frac{\Gamma(K^\perp|_P)^G}{\Gamma(K\cap
K^\perp|_P)^G}.
$$
To see that this bracket descends to a bracket on
$\Gamma(E_{red})$, one must check that if $v$ and $w$ are
$G$-invariant sections in $K^\perp|_P$ and $K^\perp\cap K|_P$,
respectively, then their bracket lies in $K^\perp\cap K|_P$. Note
that it suffices to check that their bracket lies in $K|_P$, since
we already know that it is an invariant section of $K^\perp|_P$.

 Writing the extensions of $v$ and $w$ to $\Gamma(E)$ as
$\tilde v$ and $\tilde w=\sum f_i\rho(a_i)$, we have
\begin{eqnarray*}
\Cour{\tilde v,\tilde w}&=&\sum_i (f_i\Cour{\tilde v,\rho(a_i)} +
(\pi(\tilde v)f) \rho(a_i))\\
&=&\sum_i(-f_i\Cour{\rho(a_i),\tilde v} + f_i\pi^*d\IP{\tilde v,
\rho(a_i)} + (\pi(\tilde v)f) \rho(a_i)).
\end{eqnarray*}
Restricting to $P$, we obtain $\Cour{\tilde v,\tilde w}|_P=\sum_i
(\pi(\tilde v)f)\rho(a_i)$, which is a section of $K|_P$, as
desired. The same conclusion holds for $\Cour{\tilde w,\tilde v}$,
and hence we obtain a Courant bracket on $\Gamma(E_{red})$. This
makes $E_{red}$ into a Courant algebroid over $M_{red}=P/G$ with
anchor given by the natural projection, which is clearly
surjective.

The Courant algebroid $E_{red}$ is exact if and only if the kernel
of its anchor is isotropic.  Along $P$ this can be expressed as
the condition that $\{v \in \Kperp : \pi(v) \in \pi(K)\}$ be
isotropic in $E$.  This happens if and only if $\pi(K\cap K^\perp)
= \pi(K)\cap \pi(K^\perp)$ in $TP$.  If $K$ itself was isotropic,
then $K \subset\Kperp$, and hence the condition would be
automatically satisfied.
\end{proof}



We will give explicit examples of reduced Courant algebroids using
the construction above in Section \ref{subsec:exCA}. Note that
this construction depends upon a choice of leaf $P\subset M$ of
$\Delta_b$. We end this subsection with a discussion of the
integrability of the distributions $\Delta_s$ and $\Delta_b$.


Suppose that an extended action $\rho$ is such that $G$ acts
freely and properly on the \textit{entire manifold $M$} and
$\rho(\Hh)$ has constant rank \textit{everywhere in $M$}. Then
$\Delta_s$ has constant rank by \eqref{kerhh}, and its
integrability follows from the fact that the space of sections of
$\rho(\Hh)$ is generated by closed 1-forms. Hence $\Delta_s$
defines a regular foliation of $M$. The next proposition asserts
that, although $\Delta_b$ may not have constant rank, it is a
generalized integrable distribution in the sense of Sussmann
\cite{Sus73}, defining a singular foliation of $M$ (i.e., its
leaves are smooth immersed submanifolds of varying dimensions):

\begin{prop}\label{lem:consrank}
Let $\rho:\Aa\to \Gamma(E)$ be an extended $G$-action on an exact
Courant algebroid $E$ over $M$. Assume that the $G$-action on $M$ is
free and proper and that $\rho(\Hh)$ has constant rank everywhere in
$M$. Let $S$ be a leaf of $\Delta_s$. Then the distribution $TS\cap
\psi(\Gg)$ has constant rank along $S$, and $\Delta_b$ is an
integrable generalized distribution.
\end{prop}

\begin{proof}

Recall from the proof of Lemma~\ref{cstrk} that
$\IP{\rho(a),\rho(b)}$ is constant along $S$, so $(a,b)\mapsto
\IP{\rho(a),\rho(b)}|_S$ defines a symmetric bilinear form on
$\Aa$. We consider the subspace $\Hh^\perp \subset \Aa$, and its
projection to $\Gg$, $\pi(\Hh^\perp)\subset \Gg$. We claim that
$\psi(\pi(\Hh^\perp))=TS\cap \psi(\Gg)$. Indeed, on the one hand,
$$
\psi(\pi(\Hh^\perp))=\{\psi(g)\;|\;\exists\; a\in \Aa \mbox{ with
} g=\pi(a), \mbox{ and } \rho(a)\in \rho(\Hh)^\perp\}.
$$
On the other hand, $TS\cap \psi(\Gg)=\{\psi(g)\;|\; \exists \; v
\in K^\perp, \;\psi(g)=\pi(v)\}$. But $g=\pi(a)$ for some $a\in
\Aa$, and $\psi(g)=\pi(\rho(a))$. It follows that $v-\rho(a)\in
T^*M$, i.e., $\rho(a)\in K^\perp+T^*M=\rho(\Hh)^\perp$. Hence
$TS\cap \psi(\Gg)$ has constant rank along $S$.

Now let $q:M\to M/G$ be the quotient map, which is a surjective
submersion.  Then $q|_S:S\to M/G$ has constant rank, so $q(S)$ is
an (immersed) submanifold of $M/G$. Then $P=q^{-1}(q(S))$, the
$G$-orbit of $S$, is an (immersed) submanifold of $M$ whose
tangent bundle is $\Delta_b|_P$. Hence $\Delta_b$ is integrable is
Susmann's sense.
\end{proof}
As a corollary, observe that the reduced manifolds $q(S)=q(P)=P/G$
are integral submanifolds of the smooth generalized distribution
\begin{equation}\label{redfoli}
\left.\frac{\pi(K+K^\perp)}{\pi(K)}\right/G\ \subset T(M/G),
\end{equation}
showing that this is also integrable in Susmann's sense. Therefore
the orbit space $M/G$ admits a singular foliation by submanifolds
which support the reduced Courant algebroids.

Note, however, that we do not need the global integrability of
$\Delta_s$ and $\Delta_b$ for the general construction of reduced
Courant algebroids.

\subsection{Isotropy action and the  reduced \v{S}evera class}\label{subsec:severaclass}

In this section we will give an alternative construction of
reduced Courant algebroids which clarifies condition \eqref{ex}
and allows us to describe the \v{S}evera class of an exact reduced
Courant algebroid. We start by considering the important special
case of a trivially extended action, in the sense of
Def.~\ref{def:trivially}. As the next example shows, in this case
condition~\eqref{ex} is precisely the requirement that the action
is \emph{isotropic}, i.e. $K\subset K^\perp$ or equivalently, in the language of Theorem \ref{eqphic}, the symmetric form $c$ vanishes and $\Phi = H +\xi_a$ is equivariantly closed.

\begin{example}\label{triviso}
Let $\rho:\Gg\into \Gamma(E)$ be a trivially extended action of a
free and proper action of $G$ on the manifold $M$, so that
$\Hh=\{0\}$. Then by equation~\eqref{kerhh}, we obtain
$\pi(K^\perp)=TM$, and in particular, $\Delta_s=\Delta_b=TM$.
Hence by Theorem~\ref{courred}, we obtain an exact reduced Courant
algebroid $E_{red}$ over $M_{red}=M/G$ if and only if $\pi(K) =
\pi(K\cap K^\perp)$, which occurs if and only if $K\subset
K^\perp$, since $K\cap T^*M=\{0\}$.  This provides an alternate
motivation for the requirement in~\cite{HullSpence} that $K$ be
isotropic.
\end{example}

In the case of a trivially extended, isotropic action of a compact
Lie group, we obtain the following description of the \v{S}evera
class of the reduced Courant algebroid, which appeared implicitly
in~\cite{FFS} in the context of gauging the Wess-Zumino term:
\begin{prop}\label{prop:reducedclass}
Let $G$ be a compact Lie group acting freely and properly on $M$,
and $\rho$ be a trivially extended, isotropic action on the exact
Courant algebroid $E$ over $M$.  Then if $[H]\in H^3(M,\R)$ is the
\v{S}evera class of $E$, the reduced Courant algebroid has
\v{S}evera class $q_*[\Phi]$, where $\Phi=H +\xi_a$ is the closed
equivariant extension induced by $\rho$, and $q_*$ is the natural
isomorphism
\begin{equation*}
\xymatrix{H^3_G(M,\R)\ar[r]^{q_*}& H^3(B,\R)}.
\end{equation*}
Furthermore, a splitting $\nabla:TM\into E$ induces a splitting of
$E_{red}$ if and only if $\rho(\Gg)\subset \nabla(TM)$.
\end{prop}

\begin{proof}
Since the action is isotropic, the reduced Courant algebroid is
exact, so it fits into the exact sequence
\begin{equation}\label{eq:exact}
T^*B \stackrel{dq^*}{\longrightarrow} E_{red}=(K^\perp/K)/G
\stackrel{dq}{\longrightarrow} TB,
\end{equation}
where $q:M\to B$ is the quotient map. We can find the reduced
\v{S}evera class by choosing an isotropic splitting of this
sequence. To find such a splitting, let us first choose a
$G$-invariant splitting of $E$, so that $E=(TM\oplus
T^*M,\IP{\cdot,\cdot},[\cdot,\cdot]_H)$, with $\rho(a) = X_{a} +
\xi_a$ for $a\in\Gg$,  and $i_{X_a}H = d\xi_{a}$, as in
Equation~\eqref{equi}.  Now let $\theta\in\Omega^1(M,\Gg)$ be a
connection for the principal $G$-bundle $M$. The image of the
natural map $TB \to TM/G$, $Y\mapsto Y^h$, where $Y^h$ is the
horizontal lift of $Y$, may not lie in $K^\perp/G$, but this holds
for the map
\begin{equation}\label{eq:lift}
Y\mapsto Y^h+ i_{Y^h}\IP{\theta\wedge\xi},
\end{equation}
where the 2-form $\IP{\theta\wedge\xi}$ is obtained by wedging
$\theta$ with $\xi\in\Omega^1(M,\Gg^*)$ and taking the trace.
(This 2-form is invariant since $\mathcal{L}_{X_a}\theta =
-\mathrm{ad}^*_a\theta$ and, by Equation~\eqref{actex},
$\mathcal{L}_{X_a}\xi = \mathrm{ad}_{a}\xi$.) Also note that the
image of \eqref{eq:lift} is isotropic in $K^\perp$ and intersects
$K$ trivially. Therefore \eqref{eq:lift} induces a map
$$
\nabla:TB\to E_{red}
$$
which is an isotropic splitting of \eqref{eq:exact}. Now the
induced 3-form on $TB\oplus T^*B$ is given by
\begin{align*}
\tilde{H}(X,Y,Z) & =2\IP{[\nabla(X),\nabla(Y)],\nabla(Z)}\\
&=2\IP{\Cour{X^h+i_{X^h}\IP{\theta\wedge\xi},Y^h+i_{Y^h}\IP{\theta\wedge\xi}}_H, Z^h+i_{Z^h}\IP{\theta\wedge\xi}}\\
&=2\IP{\Cour{X^h,Y^h}_{H+d\IP{\theta\wedge\xi}},Z^h}\\
&=(h+\IP{F\wedge\xi})(X,Y,Z),
\end{align*}
where in the last equality, $h$ is the basic component of $H$,
$F\in\Omega^2(M,\Gg)$ is the curvature of $\theta$, and we have
used the fact that, when evaluated on horizontal vectors,
\[
d\IP{\theta\wedge\xi}(X^h,Y^h,Z^h) = \IP{F\wedge\xi}(X^h,Y^h,Z^h).
\]
The mapping obtained here, which sends $H+\xi_a$ to the closed
form $h+\IP{F\wedge\xi}\in \Omega^3(B,\R)$ on the base, is exactly
the form-level pushdown isomorphism in equivariant cohomology:
\begin{equation*}
\xymatrix{H^3_{G}(M,\R)\ar[r]^{q_*}& H^3(B,\R)}.
\end{equation*}
So the curvature of the reduced exact Courant algebroid is
precisely the pushforward of the equivariant extension of the
original curvature induced by the extended action.

Note also that the splitting of $E_{red}$ used to calculate
$\tilde H$ depends on the choice of connection unless $\xi=0$, in
which case it is naturally induced from the original splitting of
$E$.
\end{proof}

We now explain how one can use these results about trivially
extended actions to tackle the general case. The key observation
is that there is an alternate construction of reduced Courant
algebroids which consists of two steps: first a restriction to a
small leaf $S\subset M$, and then a reduction through a trivially
extended action of a smaller group that acts on $S$.

The first step is based on the fact exact Courant algebroids may
always be \emph{pulled back} to submanifolds:

\begin{lemma}\label{lem:pullback}
Let $\iota:S\hookrightarrow M$ be a submanifold of a manifold
equipped with an exact Courant algebroid $E$. Then the vector
bundle
\begin{equation}\label{pulbak}
E_S := \frac{(\Ann(TS))^\perp}{\Ann
(TS)}=\frac{\pi^{-1}(TS)}{\Ann(TS)}
\end{equation}
inherits the structure of an exact Courant algebroid over $S$ with
\v{S}evera class $\iota^*[H]$, where $[H]$ is the class of $E$.
\end{lemma}

\begin{proof}
The subbundle $\Ann(TS)\subset T^*M\subset E$ is isotropic, so
$E_S$ has a natural nondegenerate symmetric bilinear form. It
inherits a Courant bracket by restriction, as in the proof of
Theorem~\ref{courred}, and a simple dimension count shows that
this Courant algebroid must be exact.

If a splitting $TM\into E$ were chosen, rendering $E$ isomorphic
to $(TM\oplus T^*M,\IP{\cdot,\cdot},\Cour{\cdot,\cdot}_H)$ with
$H\in \Omega^3_{cl}(M)$, then $\pi^{-1}(TS)=TS + T^*M$, and we
would obtain a natural splitting
\begin{equation*}
E_{red} = TS\oplus T^*M/\Ann(TS) = TS\oplus T^*S.
\end{equation*}
With this identification, the 3-form twisting the Courant
algebroid structure on $TS \oplus T^*S$ is simply the pull-back
$\iota^*H$.
\end{proof}

Let us consider an extended action on an exact Courant algebroid
$E$ over $M$, and let $S$ be a leaf of the small distribution
$\Delta_s$. As we saw from~\eqref{cab}, $\IP{\rho(a),\rho(b)}$ is
constant along a small leaf $S$ and induces a symmetric bilinear
form on the Courant algebra $\Aa$, for which $\Hh$ is isotropic.
Therefore we may define $\Aa_s=\Hh^\perp$ and $\Gg_s =
\pi(\Aa_s)$, noting that $\Aa_s$ is closed under the Courant
bracket.  This implies that $\Gg_s$ is a Lie subalgebra of $\Gg$,
which we call the \emph{isotropy subalgebra}, and it inherits a
symmetric bilinear form $c_s\in S^2(\Gg_s^*)$ by construction.
Therefore we obtain the sub-Courant algebra
\begin{equation*}
\xymatrix{0\ar[r]&\Hh\ar[r]&\Aa_s\ar[r]^\pi&\Gg_s\ar[r]&0},
\end{equation*}
which is mapped via the extended action $\rho$ into
$\pi^{-1}(TS)$. Quotienting by $\Hh$, we obtain a trivially
extended action $\rho_s$ of the isotropy subalgebra on the
pullback Courant algebra $E_S$ over $S$,
\begin{equation*}
\xymatrix{&0\ar[r]&\Gg_s\ar[r]^\pi\ar[d]^{\rho_s}&\Gg_s\ar[d]^\psi\ar[r]&0\\
0\ar[r]&\Gamma(T^*S)\ar[r]&\Gamma(E_S)\ar[r]&\Gamma(TS)\ar[r]&0}
\end{equation*}
which satisfies $\IP{\rho_s(a),\rho_s(b)}=c_s(a,b)$ by
construction. Note that the underlying group action on $S$ is by
the subgroup $G_s\subset G$ stabilizing $S$, which we call the
\emph{isotropy subgroup}.  Also there is a natural isomorphism
$S/G_s\into P/G$ if $P$ is a leaf of $\Delta_b$ containing $S$ and
satisfying the conditions of Theorem~\ref{courred} (see Lemma
\ref{lem:consrank}).

These arguments show that after pullback to $S$, we obtain a
trivially extended action as in Example~\ref{triviso}.  The
quotient of this pullback turns out to be naturally isomorphic to
the quotient Courant algebroid $E_{red}$ constructed in
Theorem~\ref{courred}, and we conclude that $E_{red}$ is exact if
and only if the action $\rho_s$ is isotropic, i.e.
$\rho_s(\Gg_s)\subset\rho_s(\Gg_s)^\perp$.
\begin{prop}\label{isotropk}
Let $P$ be as in Theorem~\ref{courred}, and let
$\iota:S\hookrightarrow P$ be a leaf of $\Delta_s$.  Then the
reduced Courant algebroid $E_{red}$ over $P/G$ is naturally
isomorphic to the quotient of the pullback $E_S$ by the isotropy
action $\rho_s$.  In particular, $E_{red}$ is exact if and only if
$\rho_s$ is isotropic, i.e. $c_s\in S^2(\Gg_s^*)$ vanishes.
\end{prop}
\begin{proof}The image of the isotropy action $\rho_s$ in $E_S$ is given by
\begin{equation*}
K_s =  \left.\frac{K\cap(K^\perp+T^*M)}{K\cap T^*M} \right |_S
\subset E_S= \left. \frac{K^\perp + T^*M}{K\cap T^*M}\right |_S.
\end{equation*}
Then the reduced Courant algebroid over $S/G_s$ is the $G_s$
quotient of the bundle
\begin{equation*}
\left. \frac{K^\perp_s}{K_s\cap K^\perp_s} = \frac{(K^\perp +
K\cap T^*M)/K\cap T^*M}{(K\cap K^\perp + K\cap T^*M)/K\cap
T^*M}\right |_S,
\end{equation*}
which is canonically isomorphic to $E_{red}=\left.\big
(\frac{K^\perp}{K\cap K^\perp}\big |_P\big )\right/G$ as a Courant
algebroid. Since $\rho_s$ is a trivially extended action, we
conclude from Example~\ref{triviso} that $E_{red}$ is exact if and
only if $K_s$ is isotropic in $E_S$, a condition equivalent to the
requirement that $\tilde K\subset E$ is isotropic along $P$, where
\begin{equation}\label{newk}
\tilde K = K\cap(K^\perp+T^*M).
\end{equation}
\end{proof}

The previous proposition, combined with
Prop.~\ref{prop:reducedclass} and Lemma~\ref{lem:pullback},
provides a description of the \v{S}evera class of any exact
reduced Courant algebroid: simply pull back the 3-form to the leaf
$S$ of $\Delta_s$ and apply Prop.~\ref{prop:reducedclass} for the
isotropy action $\rho_s$.

In the presence of a moment map $\mu:M\into \Hh^*$ for the
generalized action, the moment map condition $d(\mu_h)=\rho(h)$
implies that
\begin{equation*}
\ker (d\mu) = \Ann(\rho(\frak{h}))=\Delta_s,
\end{equation*}
so that the leaves of the small distribution $\Delta_s$ are
precisely the level sets $\mu^{-1}(\lambda)$ of the moment map.
Similarly the leaves of the big distribution are inverse images
$\mu^{-1}(\mathcal{O}_\lambda)$ of orbits
$\mathcal{O}_\lambda\subset \Hh^*$ of the action of $G$.  The
small leaf $S=\mu^{-1}(\lambda)$ then has isotropy Lie algebra
$\Gg_s=\Gg_\lambda$, which is the Lie algebra of $G_\lambda$, the
subgroup stabilizing $\lambda$ under the action of $G$ on $\Hh^*$.
Applying Theorem~\ref{courred} together with
Proposition~\ref{isotropk}, we obtain the following formulation of
the reduction procedure:
\begin{prop}[Moment map reduction]\label{prop:moment}
Let the extended action $\rho$ on the Courant algebroid $E$ have
moment map $\mu$. Then the reduced Courant algebroid associated to
the regular value $\lambda\in \Hh^*$ is obtained via pullback
$E_S$ along $\iota:S=\mu^{-1}(\lambda)\hookrightarrow M$, followed
by reduction by the isotropy action $\rho_\lambda$ of $G_\lambda$
on the level set, which we assume is free and proper. The result
is an exact Courant algebroid if and only if $\rho_\lambda$ is
isotropic, i.e. the induced symmetric form $c_\lambda\in
S^2(\Gg_\lambda^*)$ vanishes.
\end{prop}


\subsection{Examples}\label{subsec:exCA}

In this section we will provide some examples of Courant algebroid
reduction, illustrating the results of Sections~\ref{subsec:redCA}
and \ref{subsec:severaclass}.

\begin{ex}{1form}
Even a trivial group action may be extended by 1-forms; consider
the extended action $\rho:\RR\into \Gamma(E)$ on an exact Courant
algebroid $E$ over $M$ given by $\rho(1)=\xi$ for some closed
1-form $\xi$. Then $K = \IP{\xi}$ and $\Kperp = \{v \in \E :
\pi(v) \in \Ann(\xi)\}$ which induces the distribution
$\Delta_b=\Delta_s=\Ann(\xi)\subset TM$, which is integrable
wherever $\xi$ is nonzero. Since the group action is trivial, a
reduced Courant algebroid is simply a choice of integral
submanifold $\iota:S\hookrightarrow M$ for $\xi$ together with the
pullback exact Courant algebroid $E_{red}=E_S = K^\perp/K$, as in
Lemma~\ref{lem:pullback}. The \v{S}evera class in this case is the
pullback to $S$ of the class of $E$.
\end{ex}

\begin{ex}{vectors}
At another extreme, consider a free and proper action of $G$ on
$M$, with infinitesimal action $\psi:\Gg\into \Gamma(TM)$, and
extend it trivially by inclusion to a split Courant algebroid
$(TM\oplus T^*M,\IP{\cdot,\cdot},\Cour{\cdot,\cdot}_H)$ such that
the splitting is preserved by the action.  By
Equation~\eqref{equi}, this is equivalent to the requirement that
$H$ is an invariant basic form.

Then $K=\psi(\Gg)$ and $K^\perp = TM\oplus \Ann(K)$, so that
$\Delta_s=\Delta_b=TM$ and the reduced Courant algebroid is
\begin{equation*}
TM/K\oplus  \Ann (K) =TB\oplus T^*B,
\end{equation*}
where $B=M/G$ is the quotient and the 3-form twisting the Courant
bracket on $B$ is the push-down of the basic form $H$.
\end{ex}

The next example shows explicitly that a trivial twisting $[H]=0$
may give rise to a cohomologically nontrivial reduced Courant
algebroid.

\begin{ex}{S1 bundles}
Consider $M=S^3\times S^1$ as an $S^1$-bundle over $S^2\times
S^1$, where the $S^1$-action on the first factor of $M$ generates
the Hopf bundle $S^3\to S^2$, and the action on the second factor
is trivial. We denote the infinitesimal generator of the action on
$M$ by $\partial_t$. If $\xi$ is a volume form in $S^1$, then
$\rho(1)=\partial_t + \xi$ defines a trivially extended, isotropic
$S^1$-action on the Courant algebroid $TM\oplus T^*M$, with $H=0$.
By Prop.~\ref{prop:reducedclass}, the reduced Courant algebroid
over $S^2\times S^1$ has curvature $F\wedge \xi$, where $F$ is the
Chern class of the Hopf fibration, so the reduced \v{S}evera class
is nontrivial.
\end{ex}

As we saw in Prop.~\ref{prop:reducedclass}, a reduced Courant
algebroid may not inherit a canonical splitting. The next example
illustrates a situation where the reduced Courant algebroids are
naturally split.

\begin{example}\label{spli}
One situation where $E_{red}$ always inherits a splitting is when
$E$ is equipped with a $G$-invariant splitting $\nabla$ and the
action $\rho$ is split, in the sense that there is a splitting $s$
for $\pi:\Aa\into\Gg$ making the diagram commutative:
\begin{equation}
\xymatrix{\Aa\ar[d]^\rho&\Gg\ar[l]^s\ar[d]^\psi\\
\Gamma(E)&\Gamma(TM)\ar[l]^\nabla}
\end{equation}
In this case, the image distribution $\rho(\Aa)=K$ decomposes as
$K=K_T\oplus K_{T^*}$, with $K_T\subset TM$ and $K_{T^*}\subset
T^*M$, and hence we have the pointwise identification
\begin{equation}\label{eq:splitE}
\frac{K^\perp}{K\cap K^\perp} = \left(\frac{\Ann(K_T)}{\Ann
(K_T)\cap K_{T^*}}\right) \oplus\left( \frac{\Ann (K_{T^*})}{\Ann
(K_{T^*})\cap K_{T}}\right).
\end{equation}
Upon restriction to a leaf $P$ and quotient by $G$, the
distribution \eqref{eq:splitE} agrees with $TM_{red}^*\oplus
TM_{red}$, since $\Ann (K_{T^*})/(\Ann (K_{T^*})\cap K_T) =
\Delta_s/\rho(\Gg_s)$. Hence $E_{red}$ is split. The curvature $H$
of the given splitting for $E$ is then basic, and the curvature
for $E_{red}$ is simply the pullback to $S$ followed by pushdown
to $S/G_{s}$.

If we are in the situation above, where the action is split, one
has a natural trivially extended $G$-action on $M$ coming from
$\rho\circ s$.  Assuming that $G$ acts freely and properly on all
of $M$, we may form the quotient Courant algebroid $E_{red}^1$
over $M/G$, which is exact since $\rho\circ s$ is isotropic.
Assuming that $\rho(\Hh)$ had constant rank on $M$, then as we saw
in Section \ref{subsec:redCA}, $M/G$ inherits a generalized
foliation; the pullback of $E_{red}^1$ to a leaf of this foliation
would then recover the reduced Courant algebroid $E_{red}$ over
$M_{red}$ constructed as before.
\end{example}

An example of such a split action, where the reduced Courant
algebroid may be obtained in two equivalent ways, is the case of a
symplectic action, as introduced in Example~\ref{sym}.
\begin{ex}{symplectic action}\label{ex:sympred}
Let $(M,\omega)$ be a symplectic manifold and consider the
extended $G$-action $\rho:\Gg\oplus\Gg\into \Gamma(TM\oplus T^*M)$
with curvature $H=0$ defined in Example \ref{sym}. This is clearly
a split action in the above sense.

Let $\psi:\Gg\into \Gamma(TM)$ be the infinitesimal action and
$\psi(\Gg)^\omega$ denote the symplectic orthogonal of the image
distribution $\psi(\Gg)$.  Then the extended action has image
\begin{equation*}
K=\psi(\frak{g})\oplus \omega(\psi(\frak{g})),
\end{equation*}
so that the orthogonal complement is
\begin{equation*}
K^\perp = \psi(\Gg)^\omega \oplus \Ann(\psi(\Gg)).
\end{equation*}
Then the big and small distributions on $M$ are
\begin{equation*}
\begin{split}
\Delta_s &= \psi(\Gg)^\omega,\\
\Delta_b &= \psi(\Gg)^\omega + \psi(\Gg).
\end{split}\end{equation*}
If the action is Hamiltonian, with moment map $\mu:M\into \Gg^*$,
then $\Delta_s$ is the tangent distribution to the level sets
$\mu^{-1}(\lambda)$ while $\Delta_b$ is the tangent distribution
to the sets $\mu^{-1}(\mathcal{O}_\lambda)$, for
$\mathcal{O}_\lambda$ a coadjoint orbit containing $\lambda$.
Therefore we see that the reduced Courant algebroid is simply
$TM_{red}\oplus T^*M_{red}$ with $H=0$, for the usual symplectic
reduced space
$M_{red}=\mu^{-1}(\mathcal{O}_\lambda)/G=\mu^{-1}(\lambda)/G_{\lambda}$.

Since the action is split, we may also observe, assuming that $G$
acts freely and properly on $M$, that the quotient $M/G$ is
foliated via~\eqref{redfoli} by the possible reduced spaces.  This
generalized distribution is given in this case by
\begin{equation*}
\left.\frac{\psi(\Gg)^\omega + \psi(\Gg)}{\psi(\Gg)}\right/G\ \
=dq(\psi(\frak{g})^\omega)  \subset T(M/G),
\end{equation*}
where $q:M\to M/G$ is the quotient map. This is precisely the
distribution defined by the image of the Poisson tensor
$\Pi:T^*(M/G)\to T(M/G)$ induced by $\omega$ (recall that
$\Pi(df)=dq(X_{q^*f})$, where $X_{q^*f}$ is the Hamiltonian vector
field for $q^*f$). So a reduced manifold for the extended action
is just a symplectic leaf of $M/G$.
\end{ex}

Finally, we present an example of a reduced Courant algebroid
which is not exact.

\begin{ex}{nonexact}
Let  $\rho:\mathfrak{s}^1 \into \E$ be a trivially extended $S^1$
action which is not isotropic, i.e. $\IP{\rho(1),\rho(1)}\neq 0$.
Hence the reduced manifold for this action is just $M/S^1$ and the
reduced algebroid is $\E_{red} = (K^{\perp}/(K\cap\Kperp))/S^1$.
However, $K\cap\Kperp = \{0\}$ and so $\E_{red}$ is odd
dimensional; hence it is not an exact Courant algebroid.
\end{ex}


\section{Reduction of Dirac structures}\label{four}

A \emph{Dirac structure} \cite{Courant,LWX97} on a manifold $M$
equipped with exact Courant algebroid $E$ is a maximal isotropic
subbundle $D\subset E$ whose sections are closed under the Courant
bracket. This last requirement is referred to as the {\it
integrability} condition for $D$. When the Courant algebroid is
split, with curvature $H\in \Omega^3_{cl}(M)$, these are usually
referred to as {\it H-twisted} Dirac structures \cite{SW01}.

For $H=0$, examples of Dirac structures on $M$ include closed
2-forms and Poisson bivector fields (in these cases $D$ is simply
the graph of the defining tensor, viewed either as a map
$\omega:TM\to T^*M$ or $\Pi: T^*M\to TM$) as well as involutive regular
distributions $F\subset TM$, in which case $D=F\oplus \Ann(F)$.

In the presence of an extended action of a connected Lie group $G$
on the Courant algebroid $E$, one may consider Dirac structures
which are $G$-invariant subbundles of $E$, a condition equivalent
to the following.
\begin{definition} A Dirac structure $D\subset E$ is preserved by an
extended action $\rho$ if and only if $[\rho(\frak{a}),
\Gamma(D)]\subset \Gamma(D)$.
\end{definition}

In this section we explain how a Dirac structure which is preserved
by an extended action may be transported from a Courant algebroid
$E$ to its reduction $E_{red}$.

\subsection{Reduction procedure}

To see how Dirac structures are transported under Courant reduction,
we first explain the map at the level of linear algebra.  Suppose
that $E$ is a real vector space equipped with a nondegenerate
symmetric bilinear form of split signature, and suppose that an
isotropic subspace $K\subset E$ is given.  Then the Courant
reduction along $K$ is defined to be $E_{red} = K^\perp/K$.
Furthermore, there is a \emph{canonical relation} between $E$ and
$E_{red}$, i.e. a maximal isotropic subspace
\begin{equation*}
\varphi_K=\left\{(x,[x])\in \overline{E}\times E_{red}\ :\ x\in
K^\perp\right\},
\end{equation*}
where $\overline{E}$ denotes $E$ with negative symmetric form. If
$D\subset E$ is a Dirac structure (i.e. a maximal isotropic
subspace), view it as a relation $D\subset\{0\}\times E$.  Then
composition (as a relation) with $\varphi_K$ defines a Dirac
structure in $E_{red}$, given by
\[
D_{red}:=\varphi_K\circ D = \frac{D\cap K^\perp+K}{K}\subset
E_{red}.
\]
In this way we obtain a reduction map on Dirac structures.  This is
entirely analogous to the reduction of Lagrangian subspaces under a
symplectic reduction, as described by Weinstein~\cite{We79} using
canonical relations in the symplectic category.

Now if $E$ is an exact Courant algebroid and $K=\rho(\Aa)$ is the
image of an extended action which is isotropic along a big leaf
$P\subset M$, then any $G$-invariant Dirac structure $D$ along $P$
gives rise to the following \emph{reduced Dirac structure}, assuming
the result is a smooth bundle:
\begin{equation*}
D_{red} = \left.\frac{(D\cap K^\perp+K)|_P}{K|_P}\right/ G\
\subset E_{red}.
\end{equation*}
Note that $D_{red}$ is smooth if $D\cap K^\perp$ (or equivalently
$D\cap K$) has constant rank over $P$.  For the proof that $D_{red}$
is integrable, see Theorem~\ref{quotient dirac structure}.

If the extended action is not isotropic, the procedure just
described must be modified.  In this case we use the result of
Proposition~\ref{isotropk} that the exact reduced Courant algebroid
$E_{red}$ can be constructed by first pulling $E$ back to a leaf
$S\subset M$ of $\Delta_s$ and then taking the quotient by the
isotropic action $\rho_s$. Over the leaf $S$, the isotropic
subbundle $\rho(\Hh)=K\cap T^*M$ determines a map of Dirac
structures from $E|_S$ to the pullback Courant algebroid $E_S$. This
is a generalization of the pullback of Dirac structures defined
in~\cite{BurRadko}. After pullback, the isotropy action
$\rho_s(\Gg_s)\subset E_S$ determines a map of Dirac structures from
$E_S$ to $E_{red}$. This is a generalization of the Dirac
pushforward~\cite{BurRadko}. The composition of these maps takes any
$G$-invariant Dirac structure $D$ along $P$ to
\begin{equation}\label{gn}
D_{red} =\left.\frac{(D\cap \widetilde K^\perp+\widetilde
K)|_P}{\widetilde K|_P}\right/ G\ \subset E_{red}.
\end{equation}
where $\widetilde K = K\cap (K^\perp+T^*M)\subset E$, as defined
in~\eqref{newk}. As a result, the reduced Dirac structure is
obtained by the same procedure as in the isotropic case, applied
to $\widetilde K$ instead of $K$.  We now show that $D_{red}$,
when smooth as a vector bundle, is automatically integrable.

\begin{theo}\label{quotient dirac structure}
Let $E$, $\rho$, and $P$ be as in Theorem~\ref{courred}, and such
that $E_{red}$ is exact over $M_{red}=P/G$.  Let the action $\rho$
preserve a Dirac structure $D \subset \E$.  Then if $D_{red}$, as
described above~\eqref{gn}, is a smooth subbundle (e.g. if
$D\cap\widetilde K$ has constant rank), it defines a Dirac structure
on the reduction $M_{red}$.
\end{theo}
\begin{proof}
The only property of $D_{red}$ that remains to be checked is
integrability. To do so, we first observe that the Courant bracket
on $E_{red} = (\widetilde{K}^\perp/\widetilde{K})|_P/G$ admits the
following description, equivalent to the one given in
Theorem~\ref{courred}. Given sections $v_1, v_2$ of $E_{red}$, let
us consider representatives in $\Gamma(\widetilde{K}^\perp|_P)^G$,
still denoted by $v_1, v_2$. Then extend them to sections $\tilde
v_1, \tilde v_2$ of $E$ over $M$, and define $\Cour{v_1,v_2}$ as
$\Cour{\tilde v_1,\tilde v_2}|_P$. Similarly to
Theorem~\ref{courred}, one can show that $\Cour{\tilde v_1,\tilde
v_2}|_P \in \Gamma(\widetilde K^{\perp}|_P)^G$, and that different
choices of extensions change the bracket by invariant sections of
$\widetilde{K}$ over $P$. Also, the bracket between elements in
$\Gamma(\widetilde K|_P)^G$ and $\Gamma(\widetilde K^\perp|_P)^G$
remains in $\Gamma(\widetilde K|_P)^G$, so there is an induced
bracket on $E_{red}$. This bracket agrees with the one defined in
Theorem~\ref{courred}.

Let $v_1,v_2 \in \Gamma((D\cap \widetilde K^\perp + \widetilde
K)|_P)^G$, thought of as representing sections of $D_{red}$. We
note that, around points of $P$ where $D\cap \widetilde
K^\perp|_P$ has locally constant rank, we can write $v_i = v_i'
+v_i''$, where $v_i'$ is an invariant local section of $D \cap
\widetilde K^{\perp}|_P$, and $v_i''$ is an invariant local
section of $\widetilde K|_P$. Then the bracket of $v_1, v_2$ is
$$
\Cour{v_1' + v_1'',v_2' + v_2''} = \Cour{v_1',v_2'} +
\Cour{v_1'',v_2'} + \Cour{v_1',v_2''} + \Cour{v_1'',v_2''}.
$$
Note that the last three terms on the right-hand side are in
$\Gamma(\tilde{K}|_P)^G$. As for the first term, we know that
it lies in $\widetilde{K}^{\perp}|_P$. But since $D$ is a vector
bundle over $M$, we can locally extend $v_i'$ to sections of $D$
away of $P$ and, using these extensions to compute the bracket, we
see that $\Cour{v_1',v_2'} \in \Gamma(D|_P)$, since $D$ is
closed under the bracket. As a result, we conclude that
$\Cour{v_1,v_2}$ is in $(D\cap\widetilde{K}^\perp + \widetilde
K)|_P$ around points where $D\cap \widetilde K^\perp|_P$ is
locally a bundle.

Since the points of $P$ where $D \cap \widetilde K^{\perp}|_P$ has
locally constant rank is an open dense set, the argument above
shows that for $v_1,v_2 \in \Gamma(D_{red})$, $\Cour{v_1,v_2}$
lies in $D_{red}$ over all points in an open dense subset of
$P/G$. But since $D_{red}$ is smooth, this implies that
$\Cour{v_1,v_2} \in \Gamma(D_{red})$, hence $D_{red}$ is
integrable.
\end{proof}

The reduction of Dirac structures works in the same way for
\textit{complex} Dirac structures, provided one replaces $K$ by
its complexification $\Kc=K\otimes \C$.

%


\section{Reduction of \gcss}\label{five}

A {\it \gcs} \cite{Gua03,Hit03} on a manifold $M$ equipped with
exact Courant algebroid $E$ is a complex structure on the vector
bundle $E$ which is orthogonal \wrt\ the bilinear pairing and
whose $+i$-eigenbundle is closed under the bracket. If the Courant
algebroid is split, with curvature $H\in \Omega_{cl}^3(M)$, a
\gcs\ on $E$ is called an {\it $H$-twisted} \gcs\ on $M$.

Since a \gcs\ is orthogonal, its $+i$-eigenbundle $L\subset E\tensor
\C = E_\C$ is a maximal isotropic subbundle. Therefore a generalized
complex structure on $E$ is equivalent to a complex Dirac structure
$L$ satisfying
\begin{equation}\label{eq:gcscond}
L \cap \overline{L} = \{0\}.
\end{equation}

The {\it type} of a \gcs\ at a point $p \in M$ is the complex
dimension of the kernel of the projection $\pi:L \into T_\C M$ at
$p$. Two basic examples of generalized complex structures on a
manifold $M$ (with $H=0$) arise as follows:

\begin{itemize}

\item Let $I:TM\to TM$ be a complex structure on $M$. Then it induces a \gcs\  on $M$
by
$$\J_I = \begin{pmatrix}
-I & 0\\
0 & I^*
\end{pmatrix}.
$$
The associated Dirac structure is $L=TM^{0,1}\oplus T^*M^{1,0}$, which has type $n$.

\item Let $\omega:TM\to T^*M$ be a symplectic structure. The induced \gcs \ is
$$\J_{\w} = \begin{pmatrix}
0 & -\w^{-1}\\ \w & 0
\end{pmatrix}.
$$
The associated Dirac structure is $L=\{X-i\omega(X)\,:\, X\in T_{\mathbb{C}}M\}$, and the type is zero.
\end{itemize}
A generalized complex structure on $M^{2n}$ is of {\it complex type}
if it has type $n$ at all points,  and it is of {\it symplectic
type} if it has type zero at all points. The reader is referred to
\cite{Gua03} for more details concerning generalized complex
structures.

\subsection{Reduction procedure}

Throughout this section, $\rho: \frak{a}\to \Gamma(E)$ denotes an
extended action of a connected Lie group $G$ on an exact Courant
algebroid $E$ over a manifold $M$. Let $K=\rho(\frak{a})$, and let
$\Kc=K\otimes\C$. We fix a leaf $P \hookrightarrow M$ of the
distribution $\Delta_b$ as in Thm.~\ref{courred} and assume that
the reduced Courant algebroid $E_{red}$ over $P/G$ is exact, which
amounts to the assumption that $\tilde K = K\cap (K^\perp+T^*M)$
is isotropic along $P$.

Suppose that the extended action $\rho$ preserves a generalized
complex structure $\J$ on $\E$, i.e.,  that the associated Dirac
structure $L\subset \E_{\mathbb{C}}$ is invariant. We consider its
reduction to $E_{red}$:
\begin{equation}\label{lred}
L_{red}= \left. \frac{(L\cap \widetilde{K}_{\C}^{\perp} +
\widetilde{K}_{\C})|_P}{\widetilde{K}_{\C}|_P}\right/ G
\end{equation}
If $L_{red}$ is a smooth vector bundle, then it determines a \gcs\
on $E_{red}$ \iff\ it satisfies
$L_{red}\cap\overline{L}_{red}=\{0\}$.

\begin{lem}\label{lem:real index zero}
The distribution $L_{red}$ satisfies $L_{red} \cap \overline{L}_{red} = \{0\} $ \iff
\begin{equation}\label{neat condition}
\J \widetilde K \cap \widetilde K^{\perp} \subset \widetilde K
\mbox{ over } P.
\end{equation}
\end{lem}
\begin{proof}
It is clear from \eqref{lred} that $L_{red} \cap \overline{L}_{red}
= \{0\}$ over the reduced manifold \iff\
\begin{equation}\label{middle condition}
(L\cap \widetilde{K}_\C^{\perp} + \widetilde K_\C) \cap
(\overline{L}\cap \widetilde{K}_{\C}^{\perp} + \widetilde K_{\C})
\subset \widetilde K_{\C} \quad \mbox{ over } P.
\end{equation}
Hence, we must prove that conditions \eqref{neat condition} and
\eqref{middle condition} are equivalent.

We first prove that \eqref{neat condition} implies \eqref{middle
condition}.  Let $v \in (L\cap \widetilde K_\C^{\bot} + \widetilde
K_{\C})\cap (\overline{L} \cap \widetilde K_{\C}^\bot + \widetilde
K_{\C})$ over a given point. Without loss of generality we can
assume that $v$ is real. Since $v \in L\cap \widetilde
K_{\C}^{\bot} +  \widetilde K_{\C}$, we can find $v_L \in L \cap
\widetilde K_{\C}^\bot$ and $v_{ \widetilde K} \in \widetilde
K_{\C}$ such that $v = v_L + v_{ \widetilde K}$. Taking
conjugates, we get that $v = \BAR{v_L}  + \BAR{v_{ \widetilde
K}}$, hence $v_L -\BAR{v_L} = \BAR{v_{ \widetilde K}} - v_{
\widetilde K}$. Applying $-i \J$, we obtain
$$
v_L + \BAR{v_L} = -i\J(\BAR{v_{ \widetilde K}} - v_{ \widetilde K}).
$$
The left hand side lies in $ \widetilde K^{\perp}$ while the right
hand side lies in $\J \widetilde K$. It follows from \eqref{neat
condition} that $v_L+\BAR{v_L}\in \widetilde{K}$, hence
$v=\frac{1}{2}(v_L+\BAR{v_L}+ v_{\widetilde K}+\BAR{v_{\widetilde
K}}) \in \widetilde{K}$, as desired.

Conversely, if \eqref{neat condition}  does not hold, i.e., there is
$v \in \J  \widetilde K \cap  \widetilde K^{\perp}$ with $v \not \in  \widetilde K$,
then $v - i \J v \in L\cap  \widetilde K_{\C}^{\perp}$ and $v + i \J v \in
\BAR{L} \cap  \widetilde K_{\C}^{\perp}$. Since $v \in \J  \widetilde K$
and $\J v \in  \widetilde K$, it follows that
 $v \in  L \cap  \widetilde K_{\C}^{\perp}+ \widetilde K_{\C}$ and
$v \in \BAR{L} \cap  \widetilde K_{\C}^{\perp}+  \widetilde K_{\C}$,
showing that $(L\cap \widetilde K_{\C}^{\bot} +  \widetilde
K_{\C})\cap (\overline{L} \cap  \widetilde K_{\C}^\bot +
\widetilde K_{\C}) \not \subset  \widetilde K_{\C}$. This
concludes the proof.
\end{proof}

If the Dirac reduction of the $+i$-eigenbundle of a \gcs\ $\J$ on
$E$ defines a \gcs\ on $E_{red}$, then we denote it by $\J^{red}$.
We now present a situation where this occurs.

\begin{theo}\label{Lcap Lbar}
Let $E$, $\rho$, and $P$ be as in Theorem~\ref{courred}, and such
that $E_{red}$ is exact over $M_{red}=P/G$.  If the action preserves
a \gcs\ $\J$ on $E$ and $\J K = K$ over $P$ then $\J$ reduces to
$E_{red}$.
\end{theo}
\begin{proof}
Let us consider the isotropic distribution in the exact Courant
algebroid given by
\begin{equation}\label{eq:Lred normal}
L' =\left. \frac{(L \cap\Kc^\perp + \Kc \cap
\Kc^\perp)|_P}{(\Kc\cap \Kc^\perp)|_P}\right /G \subset
E_{red}\otimes \mathbb{C}.
\end{equation}
One can check that $L'\subset L_{red}$, so in order to show that
$L'$ and $L_{red}$ coincide, it suffices to show that $L'$ is
maximal isotropic. This is what we will check now.

Since $\J \Kperp = \Kperp$ over $P$, it follows that $K_{\C}^\perp
= L \cap K_{\C}^\perp + \overline{L}\cap K_{\C}^\perp$. Hence we
have
$$
\frac{L \cap\Kc^\perp +\overline{L} \cap\Kc^\perp +\Kc \cap
\Kc^\perp}{\Kc\cap \Kc^\perp} = \frac{\Kc^\perp}{\Kc \cap
\Kc^\perp} \;\;\;\; \mbox{ along } P.
$$
After quotienting by $G$, this implies that $L' +\overline{L'}=
E_{red}\tensor \C$, showing that $L'$ is maximal. Hence
$L_{red}=L'$. Also note that $L \cap \Kc^\perp$ has constant rank
over $P$ and, since $\Kc \cap \Kc^\perp$ is a bundle over $P$,
this implies that $L'$ as defined in \eqref{eq:Lred normal} is
smooth.

Finally, in order to conclude that $L_{red}$ induces a \gcs\, we
must check that condition \eqref{neat condition} in
Lemma~\ref{lem:real index zero} holds:
\[
\J \widetilde{K}  \cap \widetilde{K}^\perp = K\cap( \Kperp +  \J
T^*M)  \cap (\Kperp+K\cap T^*M) \subset K\cap (\Kperp+K\cap T^*M)
= \widetilde{K}.
\]
\end{proof}

\begin{cor}
If the hypotheses of the previous theorem hold and the extended
action  has a moment map $\mu:M\to \frak{h}^*$, then the reduced
Courant algebroid over $\mu^{-1}(O_{\lambda})/G$ obtains a
generalized complex structure.
\end{cor}


Theorem~\ref{Lcap Lbar} uses the compatibility condition $\J K =
K$ for the reduction of $\J$. We now observe that the reduction
procedure also works in an extreme opposite situation.

\begin{prop}\label{nondegenerate}
Let $E$, $\rho$, and $P$ be as in Theorem~\ref{courred}, and such
that $\rho(\Aa)=K$ is isotropic over $P$.  If $\rho$ preserves $J$,
and $\IP{\cdot, \cdot}$ is a nondegenerate pairing between $K$ and
$\J K$, then $\J$ reduces to $E_{red}$.
\end{prop}

\begin{proof}
As $K$ is isotropic over $P$, the reduced Courant algebroid is
exact and  $\widetilde K = K$. The nondegeneracy assumption
implies that $\J K \cap K^{\perp} = \{0\}$, and it follows that $L
\cap \Kc^\perp$ is a bundle and the Dirac reduction of $L$ is
smooth. Finally, \eqref{neat condition} holds trivially.
\end{proof}


\subsection{Symplectic structures}

We now present two examples of reduction obtained from a symplectic
manifold $(M,\w)$: First, we show that ordinary symplectic reduction
is a particular case of our construction; the second example
illustrates how one can obtain a type 1 \gcs\ as the reduction of an
ordinary symplectic structure. In both examples, the initial Courant
algebroid is just $TM \oplus T^*M$ with $H=0$.

\begin{ex}{ordinary symplectic quotient}{\em (Ordinary symplectic reduction)}
Let $(M,\w)$ be a symplectic manifold, and let $\J_{\w}$ be the
\gcs\ associated with $\omega$. Following Example~\ref{sym} and keeping
the same notation, consider a symplectic $G$-action on $M$,
regarded as an extended action. It is clear that $\J_{\omega}K=K$,
so we are in the situation of Theorem~\ref{Lcap Lbar}.

Following Example~\ref{ex:sympred}, let $S$ be a leaf of the distribution
$\Delta_s=\psi(\frak{g})^{\omega}$. Since $K$ splits as $K_T
\oplus K_{T^*}$, the reduction procedure of Theorem~\ref{quotient dirac structure} in this
case amounts to the usual pull-back of $\omega$ to $S$, followed
by a Dirac push-forward to $S/G_s=M_{red}$. If
the symplectic action admits a moment map $\mu:M\to \frak{g}^*$,
then the leaves of $\Delta_s$ are level sets $\mu^{-1}(\lambda)$,
and Theorem~\ref{Lcap Lbar} simply reproduces the usual Marsden-Weinstein
quotient $\mu^{-1}(\lambda)/G_{\lambda}$.

If the symplectic $G$-action on $M$ is free and proper, then
$\omega$ induces a Poisson structure $\Pi$ on $M/G$. We saw in
Example~\ref{ex:sympred} that the reduced manifolds fit into a
singular foliation of $M/G$, which coincides with the symplectic
foliation of $\Pi$. Following the remark at the end of Section~\ref{four}, the reduction of $\J_{\omega}$ to each leaf
can be obtained by the Dirac push-forward of $\omega$ to $M/G$,
which is just $\Pi$, followed by the Dirac pull-back of $\Pi$ to
the leaf, which is the symplectic structure induced by $\Pi$ on
that leaf.
\end{ex}

Next, we show that by allowing the projection $\pi:K \into TM$
to be injective, one can reduce a symplectic structure (type 0) to
a \gcs\ with nonzero type.

 \begin{ex}{type change up}
Assume that $X$ and $Y$ are linearly independent symplectic vector
fields generating a $T^2$-action on $M$. Assume further that
$\w(X,Y) =0$ and consider the extended $T^2$-action on $TM \oplus
T^*M$ defined by
$$
\rho(\alpha_1) = X + \w(Y);\qquad \rho(\alpha_2) = -Y + \w(X),
$$
where $\{\alpha_1,\alpha_2\}$ is the standard basis of
$\mathfrak{t}^2=\mathbb{R}^2$. It follows from $\w(X,Y) =0 $ and
the fact that the vector fields $X$ and $Y$ are symplectic that
this is an extended action with isotropic $K$.

Since $\J_{\w} K = K$, Theorem \ref{Lcap Lbar} implies that the
quotient $M/T^2$ has an induced \gcs. Note that
\[
L \cap K^\bot_\C = \{ Z - i\w(Z): Z \in \Ann(\w(X)\wedge\w(Y))   \},
\]
and it is simple to check that $X- i\w(X) \in L \cap K^\bot_\C$
represents a nonzero element in $L_{red}=(({L\cap K_\C^\bot +
K_\C})/{K_\C})/G$ which lies in the kernel of the projection
$L_{red} \to T(M/T^2)$. As a result, this reduced \gcs\ has type 1.

One can find concrete examples illustrating this construction by
considering symplectic manifolds which are $T^2$-principal bundles
with lagrangian fibres, such as $T^2\times T^2$, or the
Kodaira--Thurston manifold. In these cases, the reduced \gcs\
determines a complex structure on the base 2-torus.
\end{ex}

\subsection{Complex structures}
In this section we show how a complex manifold $(M,I)$ may have
different types of generalized complex reductions.

\begin{ex}{GIT}{\em (Holomorphic quotient)}
Let $G$ be a complex Lie group acting holomorphically on $(M,I)$,
so that the induced infinitesimal map $\rho:\Gg \into
\Gamma(TM)$ is a holomorphic map. Since $K = \rho(\Gg) \subset TM$,
it is clear that $K$ is isotropic and the reduced Courant
 algebroid is exact. Furthermore,
as $\rho$ is holomorphic, it follows that $\J_I K = K$. By Theorem
\ref{Lcap Lbar}, the complex structure descends to a \gcs\ in the
reduced manifold $M/G$. The reduced \gcs\ is nothing but the
quotient complex structure obtained from holomorphic quotient.
\end{ex}

The previous example is a particular case of a more general fact:
if $(M,I)$ is a complex manifold, then any reduction of $\J_I$ by
an extended action satisfying $\J_IK=K$ results in a generalized
complex structure of complex type. Indeed, $T^*M_{red}$ can be
identified with
$$
\left.\frac{(K^\perp\cap T^*M + K\cap K^\perp)|_P}{(K\cap
K^\perp)|_P}\right/ G \subset E_{red}
$$
and using that $\J_I(T^*M)=T^*M$, one sees that
$\J^{red}(T^*M_{red})=T^*M_{red}$, i.e., $\J^{red}$ is of complex
type. However, using Proposition~\ref{nondegenerate}, one can
produce reductions of  complex structures which are \textit{not}
of complex type.

\begin{ex}{type change down}
Consider $\C^2$ equipped with its standard holomorphic coordinates
$(z_1=x_1+iy_1, z_2=x_2+iy_2)$, and let $\rho$ be the extended
$\R^2$-action on $\C^2$ defined by
$$
\rho(\alpha_1) = \del_{x_1} + dx_2, \qquad \rho(\alpha_2) =
\del_{y_2} + dy_1,
$$
where $\{\alpha_1,\alpha_2\}$ is the standard basis for $\R^2$. Note
that $K=\rho(\R^2)$ is isotropic, so the reduced Courant algebroid
over $\C/\R^2$ is exact. Since the natural pairing between $K$ and
$\J_I K$ is nondegenerate,  Proposition~\ref{nondegenerate} implies
that one can reduce $\J_I$ by this extended action.  In this
example, one computes
$$
\Kc^\perp \cap L = \mbox{span}\{\del_{x_1} - i \del_{x_2} - dy_1
+i dx_1, \del_{y_1} - i \del_{y_2} - dy_2 + i dx_2\}
$$
and $ \Kc^\perp \cap L \cap \Kc = \{0\}$. As a result, $L_{red}
\cong \Kc^\perp \cap L$. So $\pi:L_{red} \into \C^2/\R^2$ is an
injection, and $\J^{red}$ has zero type, i.e., it is of symplectic
type.
\end{ex}


\subsection{Extended Hamiltonian actions}

In order to reduce a generalized complex structure $\J$ preserved by
an extended action, we saw in Theorem~\ref{Lcap Lbar} that a
sufficient condition is the compatibility $\J K=K$.  Natural
examples where this condition holds arise as follows: one starts
with an action generated by sections $v_i\in \Gamma(E)$, and then
enlarges it to a new extended action generated by sections
\begin{equation}\label{eq:complex}
\{v_i,\J v_j\}.
\end{equation}
Examples where this construction works are the extended actions
associated with symplectic and holomorphic actions: in the
symplectic case (see Example \ref{sym}), one starts with symplectic
vector fields $X_i$ and defines an extended action of the
hemisemidirect Courant algebra, by adding new generators
$\J_{\omega}(X_j)=\omega(X_j)$, which act as closed 1-forms; in the
holomorphic case, one starts with an action generated by $X_i$
preserving a complex structure $I$, and then forms the (trivially)
extended action of the complexified Lie algebra, generated by
$\{X_i,\J_I X_j\}$, where now $\J_I X_j=I X_j$ are new vector
fields.

The ``complexification'' \eqref{eq:complex} does not always define
an extended action, as we will see.  However, in the case of a
\emph{Hamiltonian} action we show that it does produce an example of
an extended action satisfying $\JJ K = K$.

It is familiar in the case of a complex manifold that a real vector
field $X$ preserves the complex structure $I$ if and only if its
$(1,0)$ component $X^{1,0}\in T_{1,0}M$ is a holomorphic vector
field. Therefore $IX= iX^{1,0}- iX^{0,1}$ also preserves the complex
structure. In particular, if $X$ generates an $S^1$ action then
$\{X,IX\}$ defines a holomorphic $\C^*$ action on the complex
manifold.

For generalized complex structures a similar phenomenon occurs,
except that symmetries are governed by the differential complex
$(\Omega^\bullet(L)=\Gamma(\wedge^\bullet L^*),d_L)$ associated to
the complex Lie algebroid $L$ defined by the $+i$--eigenbundle of
$\J$.

\begin{lemma}\label{lem:dL}
A real section $v\in \Gamma(E)$
preserves the generalized complex structure $\JJ$ under the adjoint
action if and only if $d_Lv^{0,1}=0$, where $v=v^{1,0}+v^{0,1}\in
L\oplus \overline{L}=E\otimes \C$ and we use the inner product to identify
$\overline{L}=L^*$.
\end{lemma}

\begin{proof}
A real section $v\in \Gamma(E)$ preserves $\J$ if and only if
$[v,\Gamma(L)]\subset \Gamma(L)$. Since $L$ is maximal isotropic, it
suffices to check that $\IP{[v^{0,1}, w_1],w_2}=0$ for all
$w_1,w_2\in \Gamma(L)$. By definition of the Lie algebroid
differential $d_L$, and using the basic properties of the Courant
bracket, we have
\begin{eqnarray*}
d_L v^{0,1}(w_1,w_2)& = & \pi(w_1)\IP{v^{0,1},w_2} -
\pi(w_2)\IP{v^{0,1},w_1} -
\IP{v^{0,1}, [w_1,w_2]}\\
 & = & 2\IP{[v^{0,1}, w_2], w_1} + \IP{v^{0,1},[w_1, w_2]} +
\pi(w_2)\IP{v^{0,1},w_1} - \pi(w_1)\IP{v^{0,1},w_2}\\
 &=& 2 \IP{[v^{0,1},w_2], w_1} - d_L v^{0,1}(w_1,w_2),
\end{eqnarray*}
so $d_L v^{0,1}(w_1,w_2)=\IP{[v^{0,1}, w_2], w_1}$, which
immediately implies the result.
\end{proof}

We obtain the following exact
sequence describing $\Gamma_\J(E)$, the space  of sections of $E$
preserving $\JJ$ under the adjoint action \cite{Gua03}:
\begin{equation*}
\xymatrix{C^{\infty}(M,\C)\ar[r]^{D}& \Gamma_\J(E) \ar[r]&H^1(L)\ar[r]&0},
\end{equation*}
where $D(f)= d_Lf + \overline{d_L f}\in \Gamma(E)$, and the final
term denotes the first Lie algebroid cohomology of $L$. The sections
of $E$ which lie in the image of $D$ are called \emph{Hamiltonian}
symmetries \cite{Gua03}, in direct analogy with the symplectic case.
Note that for $f\in C^{\infty}(M,\CC)$ we have by definition
\begin{equation*}
d_Lf = \tfrac{1}{2}(df + i \JJ df),
\end{equation*}
so that the operator $D$ may be expressed as
\begin{equation*}
Df = d(\mathrm{Re} f) - \JJ d(\mathrm{Im} f).
\end{equation*}
Also note that the projection $\pi(Df)\in \Gamma(TM)$ lies in the
projection $\pi(\JJ(T^*M))$ of the Dirac structure $\JJ(
T^*M)\subset E$ and hence is tangent to the symplectic leaves of the
Poisson structure induced by $\JJ$. This places a strong constraint
on Hamiltonian symmetries which is familiar from the situation in
Poisson geometry.
\begin{example}
In the symplectic case, a section $X+\xi\in \Gamma(TM\oplus T^*M)$
preserves $\J_{\omega}$ precisely when $X$ is a symplectic vector field and
$d\xi=0$, whereas it is Hamiltonian if and only if $X$ is
Hamiltonian in the usual sense and $\xi$ is exact.  In the complex
case, $X+\xi$ preserves $\J_I$ when $X^{1,0}$ is holomorphic and
$\delbar\xi^{0,1}=0$, whereas it is Hamiltonian if and only if $X=0$
and $\xi=\delbar f  + \del\overline{f}$ for $f\in C^{\infty}(M,\CC)$.
\end{example}

We have the following immediate consequence of Lemma~\ref{lem:dL}.

\begin{cor}
If $v\in \Gamma(E)$ preserves $\JJ$ then so does $\JJ v = iv^{1,0}-iv^{0,1}$.
\end{cor}

However, if the infinitesimal action of $v$ integrates to
an extended action on $E$, then this does not guarantee that $\JJ v$
also does, as we now show.
\begin{ex}{forms}
Let $\frak{h} = \frak{a} = \R$ be a Courant algebra over the trivial
Lie algebra $\Gg=\{0\}$ and consider an action by covectors
$\rho:\frak{a} \into T^*M\subset TM\oplus T^*M$. In order that
$\rho$ define an extended action we need $\xi = \rho(1) \in
\Omega_{cl}^1(M)$. If $M$ is endowed with a complex structure $I$,
then the complexification of $\rho$ satisfies $\rho_{\C}(i) =
I^*\xi$, which is closed only if $d^c\xi=0$.
\end{ex}
While the ``complexification'' proposed in~\eqref{eq:complex} may be
obstructed because of the fact that $\JJ v$ may not define an
extended action even if $v$ does, we now show that if the given
action is \emph{Hamiltonian}, then it is equivalent, in the sense of
Definition \ref{eqaction}, to an action which can be extended so
that $\JJ K = K$.

\begin{theo}\label{ham}
Let $\rho:\Gg\into \Gamma(E)$ be a trivially extended, isotropic,
Hamiltonian action on a generalized complex manifold, i.e.
$\rho(a)=D(f_a)$ for a $\Gg$-equivariant function $f:M\into
\Gg_\CC^*$. Then the equivalent action $\tilde\rho(a)=\rho(a) -
d(\mathrm{Re} f_a)$ may be extended to an action of the
hemisemidirect Courant algebra $\Gg\oplus\Gg$, with moment map
$\mathrm{Im}f$, and which satisfies the condition $\JJ K = K$.
\end{theo}
\begin{proof}
Since $\rho(a)=D(f_a)=d(\mathrm{Re} f_a) - \JJ d(\mathrm{Im} f_a)$,
we see that
\begin{equation*}
\JJ\tilde\rho(a) = d(\mathrm{Im} f_a),
\end{equation*}
which shows that the map $\rho':\Gg\oplus\Gg\into \Gamma(E)$ given
by
\begin{equation*}
\rho':(g,h)\longmapsto \tilde\rho(g) + d(\mathrm{Im} f_h)
\end{equation*}
defines an extended action, as we saw in Proposition~\ref{prop:will
need later}, and by construction satisfies $\JJ K = K$.
\end{proof}

Although this theorem concerns only Hamiltonian actions, which for
generalized complex structures is increasingly restrictive as the
type grows, we will use it to construct new examples of generalized
K\"ahler structures (see Section~\ref{six}). Also note that
Examples~\ref{type change up}, \ref{GIT} and \ref{type change down}
are not Hamiltonian. We remark that the actions which are
independently described by Lin and Tolman \cite{LT05}, as well as Hu
\cite{Hu05}, can be seen to be of this Hamiltonian type.

Finally, we provide a cohomological criterion which
determines if a given action is Hamiltonian.  If a trivially
extended, isotropic action $\rho:\Gg\into \Gamma(E)$ preserving
$\JJ$ is given, then we may decompose $\rho(a)=Z_a+\zeta_a\in
L\oplus \overline{L}$ for all $a\in \Gg_\C$.  Since this is a
Courant morphism, we may then define an equivariant Cartan model for
the differential complex $(\Omega^\bullet(L), d_L)$, by considering
equivariant polynomial functions $\Phi:\Gg_\C\lra
\Omega^\bullet(L)$, and equivariant derivative
\begin{equation*}
(d_{\Gg_\C}\Phi)(a)=d(\Phi(a))- i_{Z_a}\Phi(a),\ \ \forall
a\in\Gg_\C.
\end{equation*}
Since $\rho(a)$ preserves $\JJ$, we see that $\zeta_a\in
\Gg^*_\C\otimes \Omega^1(L)$ defines an equivariant closed 3-form.
Supposing that $[\zeta_a]=0$ in $H^3_{\Gg_\C}(L)$, we then have
\begin{equation*}
\zeta_a = d_{\Gg_\C}(\eps + h_a),
\end{equation*}
for $\eps\in\Omega^2(L)$ an invariant $d_L$-closed form and
$h_a\in\Gg^*_\CC\otimes\Omega^0(L)$ an equivariant function.
Supposing further that $[\eps]=0$ in the invariant cohomology
$H^2(L)^{\Gg_\C}$, then $\eps=d_L\eta$ for an invariant 1-form
$\eta$, and
\begin{equation*}
\zeta_a = d_{\Gg_\C}(h_a + i_{Z_a}\eta),
\end{equation*}
implying that $\rho(a)= D (f_a)\ \ \forall a\in\Gg$, where
$f_a=h_a+i_{z_a}\eta$.  This provides the following result.
\begin{prop}
Let $\rho$ be a trivially extended, isotropic action preserving a
generalized complex structure.  Then it is Hamiltonian if and only
if the classes $[\zeta_a]\in H^3_{\Gg_\C}(L)$ and $[\eps]\in
H^2(L)^{\Gg_\C}$, defined above, vanish.
\end{prop}


\section{Generalized K\"ahler reduction}\label{six}
A {\it \gks} \cite{Gua03} on an exact Courant algebroid $E$ is a
pair of commuting \gcss\ $\J_1$ and $\J_2$ such that
$$ \IP{\J_1\J_2v,v} >0 \qquad \mbox{ for all } v \in \E .$$
The symmetric endomorphism $\mc{G} = \J_1 \J_2$ therefore defines a
positive-definite metric on $E$, called the {\it \gk\ metric}.

\subsection{Reduction procedure}
In this section we follow the standard treatment of K\"ahler reduction~\cite{GuSt82,Kirwan}
 and extend it to the generalized setting.

\begin{theo}\label{gk quotient}{\em (Generalized K\"ahler reduction):}
Let $E$, $\rho$, and $P$ be as in Theorem~\ref{courred}, with
$\rho(\Aa)=K$ isotropic along $P$. If the action preserves a \gk\
structure $(\J_1,\J_2)$ and $\J_1K=K$ along $P$, Then $\J_1$ and
$\J_2$ reduce to a \gks\ on $E_{red}$.
\end{theo}

\begin{proof}
Since $K$ is isotropic, the reduced Courant algebroid is exact, and
by Theorem~\ref{Lcap Lbar}, $\J_1$ descends to $E_{red}$. In order
to show that $\J_2$ also descends, we will find an identification of
$E_{red}$ with a subbundle of $K^\perp$ which is invariant by
$\J_2$.

Let $K^\mc{G}$ denote the orthogonal of $K$ with respect to the
metric $\mc{G}$. Since $\J_1K=K$,
\begin{equation}\label{eq:orthog}
K^\mc{G}=(\J_2\J_1 K)^\perp = (\J_2 K)^\perp = \J_2K^\perp \;\;\;
\mbox{ over } P.
\end{equation}
Since $K\subset K^\perp$ along $P$, we have the
$\mc{G}$-orthogonal decomposition of $K^\perp$ as
$$
K^\perp=K\oplus (K^\mc{G}\cap K^\perp) \;\;\; \mbox{ over } P.
$$
It follows from \eqref{eq:orthog} that $K^\mc{G}\cap K^\perp$ is
$\J_2$-invariant. Using the natural identification
$$
K^\perp/K \cong K^\mc{G}\cap K^\perp \;\;\;\; \mbox{ over }P
$$
and the fact that $\J_2$ is $G$-invariant, we obtain after
quotienting by $G$ an induced orthogonal endomorphism
$\J_2^{red}:{E}_{red}\to {E}_{red}$ satisfying
$(\J_2^{red})^2=-1$. It remains to check that $\J_2^{red}$ is
integrable.

In order to verify integrability, we first describe the
$+i$-eigenbundle of $\J_2^{red}$. Let $L_2$ be the
$+i$-eigenbundle of $\J_2$. The $+i$-eigenbundle of $\J_2^{red}$
is the image under the natural projection $p:\Kc^\perp|_P \to
E_{red}\otimes \mathbb{C}$ of $L_2\cap (\Kc^\perp \cap
\Kc^\mc{G}).$ But since
$$
L_2\cap \Kc^\perp = \J_2(L_2\cap\Kc^\perp)= L_2\cap
\J_2(\Kc^\perp)=L_2\cap \Kc^\mc{G} \;\;\;\; \mbox{ over }P,
$$
it follows that the $+i$-eigenbundle of $\J_2^{red}$ is
$$
p(L\cap (\Kc^\perp \cap \Kc^\mc{G}))=p(L\cap
\Kc^\perp)=(L_2)_{red},
$$
i.e., the reduction of the Dirac structure $L_2$. It follows that
$(L_2)_{red}$ is a smooth and maximal isotropic subbundle of
$E_{red}\otimes \mathbb{C}$, and by Theorem \ref{quotient dirac
structure} we know that it is integrable. So $\J_2^{red}$ is
integrable.

Finally, we need to show that  $(\J_1^{red},\J_2^{red})$, where
$\J_1^{red}$ is the reduction of $L_1$, is a \gk\ pair in
${E}_{red}$. For that, we note that $K^\mc{G}\cap K^\perp$ is
$\J_1$-invariant along $P$, since $\J_1(K^\perp)=K^\perp$ and
$\J_1(K^\mc{G})=K^\mc{G}$. So $\J_1$ induces an endomorphism of
$E_{red}$, which coincides with the Dirac reduction $\J_1^{red}$
since they have the same $+i$-eigenbundle: indeed,
$$
\frac{L_1\cap \Kc^\perp + \Kc}{\Kc}=\frac{(L_1\cap \Kc) \oplus
(L_1\cap \Kc^\mc{G}\cap \Kc^\perp) + \Kc}{\Kc} \;\;\; \mbox{ over
} P,
$$
therefore, after quotienting by $G$, we get $(L_1)_{red}
=p(L_1\cap \Kc^\mc{G} \cap \Kc^\perp)$. The fact that $\J_1^{red}$
and $\J_2^{red}$ form a generalized K\"ahler pair is now a direct
consequence of the fact that the restrictions of $\J_1$ and $\J_2$
to $K^\mc{G}\cap K^\perp$ commute and their product is positive
definite.
\end{proof}
An important particular case of Theorem \ref{gk quotient} is when
the extended action admits a moment map.
\begin{cor}\label{cor:mom}
Let $(\J_1,\J_2)$ be a generalized K\"ahler structure preserved by
an extended action admitting a moment map $\mu:M\to \frak{h}^*$.
Assume that the $G$-action on $\mu^{-1}(0)$ is free and proper. If
$\J_1(K)=K$ over $\mu^{-1}(0)$, and the induced symmetric form $c_0
\in S^2\frak{g}^*$ vanishes, then $\J_1$ and $\J_2$ can be reduced
to $M_{red}$ and define a generalized K\"ahler structure.
\end{cor}
This corollary follows from the fact that if $c_0$ vanishes, then
both the isotropy action and the full action along $\mu^{-1}(0)$ are
isotropic, i.e. $K\subset K^\perp$ on the level set.  Of course
these hypotheses are all fulfilled for a complexified Hamiltonian
action as in Theorem~\ref{ham}.  We now state the particular case
when $\J_1$ is a symplectic structure since we use it in the next
section.
\begin{cor}\label{symplectic reduction}
Let $(\J_1,\J_2)$ be a \gk\ structure on $E=TM\oplus T^*M$ with
$H=0$, such that $\J_1$ is an ordinary symplectic structure. Assume
that there is a Hamiltonian action on $(M,\J_1)$, with moment map
$\mu:M\to \frak{g}^*$, and preserving $\J_2$. If the action of $G$
on $\mu^{-1}(0)$ is free and proper, then the symplectic reduced
space $M_{red}=\mu^{-1}(0)/G$ carries a \gk\ structure given by
$(\J_1^{red},\J_2^{red})$.
\end{cor}
This result was independently obtained in \cite{LT05}, where it is
used to produce many examples of generalized K\"ahler quotients.
Also, when $\J_2$ is a complex structure, then $\J_2^{red}$ is as
well, and we recover the original K\"ahler reduction of
\cite{GuSt82,Kirwan}.

\begin{ex}{symplectic cut}
(Symplectic cut): Let $(\J_1,\J_2)$ be a generalized Kahler manifold
as in Corollary \ref{symplectic reduction}. Assume that there is a
Hamiltonian $S^1$-action on $M$ preserving $\J_2$, and let $f:M
\into \R$ be its moment map. Consider $\C$ with its natural K\"ahler
structure $(\omega,I)$, and equipped with the $S^1$-action
$\theta\cdot z:=e^{i\theta}z$. Then $N = M \times \C$ has a \gk\
structure $(\J_1', \J_2')$, where $\J_1'$ is the product symplectic
structure and $\J_2' = \J_2 \times I$, and
$$
\mu : N \into \R; \qquad \mu(p,z) = f(p) + |z|^2
$$
is a moment map for the diagonal $S^1$-action on $N$. This action
preserves the \gks\ so, by Corollary \ref{symplectic reduction},
the symplectic quotient of $N$ inherits a \gks.
\end{ex}


\subsection{Examples of generalized K\"ahler structures on $\mathbb{C}P^2$}

Now we apply the results from the last section to produce new
examples of \gks\ on $\mathbb{C}P^2$ with type change. The method
consists of deforming the standard K\"ahler structure in $\C^3$ so
that the deformed structure is still preserved by the circle
action
\begin{equation}\label{circle action in C3}
e^{i\theta}: (z_1,z_2,z_3) \mapsto  (e^{i \theta} z_1,e^{i \theta} z_2, e^{i \theta}z_3).
\end{equation}
Then Corollary \ref{symplectic reduction} implies that $\C P^2$, regarded as a symplectic reduction of
$\C^3$, inherits a reduced \gks.

In the computations that follow, it will be convenient to use
differential forms to describe a \gcs\ $\J$ on a manifold $M$. So we
recall from \cite{Gua03} that $\J$  is completely determined by its
{\it canonical line bundle}, $C \subset \wedge^\bullet T^*_{\C}M$.
This bundle is defined as the Clifford annihilator of $L$, the
$+i$-eigenspace of $\J$. The fact that $L$ is a Dirac structure of
real index zero ($L \cap \overline{L} = \{0\}$) translates into
properties for $C$: if $\gf$ is a nonvanishing local section of $C$,
then
 \begin{itemize}
 \item At each point, $\gf = e^{B + i \w}\wedge  \gO$, where $B$ and $\w$ are real 2-forms and $\gO$ is a decomposable complex $k$-form;
 \item There is a local section $X + \xi \in \Gamma(TM \oplus T^*M)$ such that
 $$d\gf = (X + \xi)\cdot \gf = i_X \gf + \xi \wedge \gf;$$
 \item If $\sigma$ is the linear map which acts on $k$-forms by $\sigma(\ga) = (-1)^{\frac{k(k-1)}{2}}\ga$, then the Mukai
 pairing $(\gf,\overline{\gf})$ must be nonzero, where
$$ (\gf,\overline{\gf}) := (\gf \wedge \sigma(\overline{\gf}))_{top}.$$
The subscript $top$ indicates a projection to the volume form
component.
 \end{itemize}

We begin with the standard K\"ahler structure on
$(\C^3,\J_\w,\J_I)$, defined by the following differential forms:
\begin{align*}
\Omega &= dz_0 \wedge dz_1\wedge dz_2\\
\omega &= \tfrac{i}{2}(dz_0\wedge d\bar z_0 + dz_1\wedge d\bar z_1 +
dz_2\wedge d\bar z_2)
\end{align*}

As explained in \cite{Gua03}, it is possible to deform this
K\"ahler structure as a generalized K\"ahler structure in such a way
that $\omega$ is unchanged whereas the complex structure $\Omega$
becomes a generalized complex structure of generic type 1.  To
achieve this, we must select a deformation $\e \in
\Gamma(L_+^*\otimes L_-^*)$, where
\[
L_\pm^* = \{X\pm i\omega(X)\ :\ X\in TM^{1,0}\},
\]
which satisfies the Maurer-Cartan equation $\delbar\eps +
\tfrac{1}{2}[\eps,\eps]=0$.  Then in regions where $\eps$ does not
invalidate the open condition that $e^\eps\Omega$ be of real index
zero, $(e^\eps\Omega,e^{i\omega})$ will be a generalized K\"ahler
pair.

\begin{ex}{example 1} In this example we deform the structure in $\C^3$ so that the
reduced structure in $\C P^2$ has type change along a triple line. A
similar deformation and quotient has been considered independently
by Lin and Tolman in \cite{LT05}, where they also consider a variety
of other examples.  A generalized K\"ahler structure on $\C P^2$
with type change along a triple line was also recently constructed
by Hitchin \cite{Hit05} using a different method.

\vskip6pt
\noindent
{\bf The deformation.}
We select the decomposable element
\[
\eps = \tfrac{1}{2}z_0^2(\del_1 + \tfrac{1}{2}dz_1)\wedge(\del_2 -
\tfrac{1}{2}dz_2),
\]
whose bivector component $\tfrac{1}{2}z_0^2\del_1\wedge\del_2$ is a
quadratic holomorphic Poisson structure.  The projectivization of
this structure is a Poisson structure on $\C P^2$ vanishing to order
3 along the line $z_0=0$.
The deformed complex structure in $\C^3$ can be written explicitly (we omit
the wedge symbol):
\begin{align}
\gf =e^{\eps}dz_0dz_1dz_2&= (1+\eps)dz_0dz_1dz_2 \notag \\
&=dz_0dz_1dz_2- \tfrac{1}{2}z_0^2 dz_0 -\tfrac{1}{4}z_0^2
dz_0dz_2d\bar z_2+ \tfrac{1}{4}z_0^2 dz_0dz_1d\bar z_1 +
\tfrac{1}{8}z_0^2dz_0dz_1dz_2d\bar z_2d\bar z_1 \notag \\
&=-\tfrac{1}{2}z_0^2dz_0\exp(-\tfrac{2}{z_0^2} dz_1dz_2 +
\tfrac{1}{2}(dz_2d\bar z_2 -dz_1d\bar z_1)) \label{deformed rho}
\end{align}
Let $\zeta = -\tfrac{1}{2}z_0^2dz_0$ and $b+i\sigma =
-\tfrac{2}{z_0^2} dz_1dz_2 + \tfrac{1}{2}(dz_2d\bar z_2 -dz_1d\bar
z_1)$.  Then the pure differential form $\gf$ is of real index zero
as long as the Mukai pairing of $\gf$ with its complex conjugate
satisfies
\[
(\gf,\overline{\gf})=\sigma^2\wedge \zeta \wedge\overline\zeta \neq 0.
\]
Calculating this quantity, we obtain:
\[
\sigma^2\wedge \zeta \wedge\overline\zeta = \tfrac{1}{2}(4-|z_0|^4)
dz_0dz_1dz_2d\bar z_0d\bar z_1d\bar z_2,
\]
proving that $(\gf,e^{i\omega})$ defines a generalized K\"ahler
structure in $\C^3$ away from the cylinder $|z_0|=\sqrt 2$.

\vskip6pt
\noindent
{\bf The reduction.}
Notice that the line generated by $\gf$, and hence the generalized
complex structure it defines, is invariant by the $S^1$-action given by \eqref{circle action in C3}. Hence,
by Corollary \ref{symplectic reduction}, the symplectic reduction of $\C^3$ will have a reduced \gks\ induced by the deformed structure above. We spend the rest of this example describing this structure.
The particular reduction we wish to calculate is the
quotient of the unit sphere $\sum_i z_i\bar z_i=1$ by the
$S^1$-action give by \eqref{circle action in C3}.

We begin with the generalized complex structure $\gf$ given by equation \eqref{deformed rho}.
The induced Dirac structure on the reduced Courant algebroid may be
calculated by pulling back to the unit sphere in $\C^3$ and pushing
forward to the quotient. The latter operation on differential forms
may be expressed simply as interior product with $\del_\theta$, the generator of the circle action
$$ \del_\theta =i(z_0\del_0 - \bar z_0 \delbar_0+z_1\del_1 - \bar z_1 \delbar_1+z_2\del_2 - \bar z_2 \delbar_2),$$
and this commutes with pull-back to the sphere. So let us first take
interior product:
\begin{align*}
i_{\del_\theta}\gf &=(i_{\del_\theta}\zeta)\exp(\tfrac{-\zeta \wedge
i_{\del_\theta}(b+i\sigma)}{i_{\del_\theta}\zeta} + b+i\sigma)\\
&=-\tfrac{i}{2} z_0^3\exp(-\tfrac{dz_0}{z_0}(\tfrac{2(z_2dz_1 -
z_1dz_2)}{z_0^2} + \tfrac{z_2d\bar z_2 + \bar z_2 dz_2 - z_1d\bar
z_1 - \bar z_1 dz_1}{2})-\tfrac{2dz_1dz_2}{z_0^2}  +
\tfrac{dz_2d\bar z_2 -dz_1d\bar z_1}{2})
\end{align*}
Now we pull back to $S^5$ by imposing $1=R^2 = \sum_i z_i\bar z_i$
and obtain a homogeneous differential form after rescaling:
\[
\tilde\gf =\exp(-\tfrac{dz_0}{z_0}(\tfrac{2(z_2dz_1 -
z_1dz_2)}{z_0^2} + \tfrac{z_2d\bar z_2 + \bar z_2 dz_2 - z_1d\bar
z_1 - \bar z_1 dz_1}{2R^2})-\tfrac{2dz_1dz_2}{z_0^2}  +
\tfrac{dz_2d\bar z_2 -dz_1d\bar z_1}{2R^2}).
\]
The holomorphic Euler vector field is $\euler = \sum_i z_i \del_i$ and
$\del_\theta = i({\euler}-\bar{\euler})$.  The radial vector field is $\del_r =
{\euler}+\bar{\euler}$.  In order to be the pull-back of a form on $\C P^2$, a
differential form $\alpha$ on $\C^3$ must satisfy $\mc{L}_{\euler}\alpha =
\mc{L}_{\bar{\euler}}\alpha = i_{\euler}\alpha= i_{\bar{\euler}}\alpha = 0$.  We have
already ensured that $\mc{L}_{\euler}\tilde\gf = \mc{L}_{\bar {\euler}}\tilde\gf=0$ and
$i_{{\euler}-\bar{\euler}}\tilde\gf = 0$, so now we may add a multiple of $dR$
to ensure $i_{{\euler}+\bar{\euler}}\tilde\gf =0$.  Since $dR$ vanishes on the
sphere, this is a trivial modification.

Recall that $i_{{\euler}+\bar{\euler}} \tfrac{dR}{R}= 1$, so we shall subtract
\begin{align*}
\tfrac{dR}{R}\wedge i_{{\euler}+\bar{\euler}}\tilde\gf &=
\tfrac{dR}{R}(\tfrac{dz_0}{z_0}(\tfrac{z_2\bar z_2 - z_1\bar
z_1}{R^2}) + \tfrac{\bar z_1 dz_1 - \bar z_2 dz_2}{R^2})\tilde\gf
\end{align*}
Finally we get a manifestly projective representative for the
generator of the canonical bundle:

\begin{align*}
\gf_B=\exp(-\tfrac{dz_0}{z_0}(\tfrac{2(z_2dz_1 - z_1dz_2)}{z_0^2} +
\tfrac{z_2d\bar z_2 + \bar z_2 dz_2 - z_1d\bar z_1 - \bar z_1
dz_1}{2R^2})&-\tfrac{2dz_1dz_2}{z_0^2}  + \tfrac{dz_2d\bar z_2
-dz_1d\bar z_1}{2R^2}\\
&-\tfrac{dR}{R}(\tfrac{dz_0}{z_0}(\tfrac{z_2\bar z_2 - z_1\bar
z_1}{R^2}) + \tfrac{\bar z_1 dz_1 - \bar z_2 dz_2}{R^2}))\notag
\end{align*}

This differential form is closed, but blows up along the type change
locus, where one can see by rescaling that it defines a complex
structure.  This generalized complex structure, together with the
Fubini-Study symplectic structure, forms a generalized K\"ahler
structure on $\C P^2$.

It may be of interest to express this generalized K\"ahler structure
in affine coordinates $(z_1,z_2)$ where $z_0=1$.  Then the type
change locus is the line at infinity.  Define $r^2 = z_1\bar z_1 +
z_2\bar z_2$:
\[
\gf_B =\exp(-2dz_1dz_2  + \tfrac{dz_2d\bar z_2 -dz_1d\bar
z_1}{2(1+r^2)} -\tfrac{1}{2}\tfrac{d(r^2)(\bar z_1 dz_1 - \bar z_2
dz_2)}{(1+r^2)^2}))
\]
The form defining the Fubini-Study symplectic form in these
coordinates is, as usual:
\[
\gf_{A} = \exp(-\tfrac{1}{2}\tfrac{(1+r^2)(dz_1d\bar z_1 +
dz_2d\bar z_2) -(\bar z_1dz_1+ \bar z_2dz_2)(z_1d\bar z_1 + z_2d\bar
z_2))}{(1+r^2)^2})
\]
An important constituent of a generalized K\"ahler structure is its
associated bi-Hermitian metric; this can be derived from the above
forms as follows.  Define real 2-forms $\omega_1,\omega_2,b$ such
that $\gf_A =e^{i\omega_1}$ and $\gf_B=e^{b+i\omega_2}$.  Then the
bi-Hermitian metric $g$ is simply
\[
g = -\omega_2 b^{-1} \omega_1.
\]
\end{ex}

\begin{ex}{Example 2}
To demonstrate the versatility of the quotient construction we now
construct a generalized K\"ahler structure on $\C P^2$ with type
change along a slightly more general cubic:  the union of three
distinct lines forming a triangle.  We postpone the discussion of
the general cubic curve to a future paper.

\vskip6pt
\noindent
{\bf The deformation.}
In this example we select a deformation $\eps$ given by the
following decomposable section of $L_+^*\otimes L_-^*$:
\[
\eps = \tfrac{1}{2}(z_0(\del_1+\tfrac{1}{2}d\bar z_1) +
z_1(\del_2+\tfrac{1}{2}d\bar z_2)+z_2(\del_0+\tfrac{1}{2}d\bar
z_0))\wedge(z_0(\del_2-\tfrac{1}{2}d\bar
z_2)+z_1(\del_0-\tfrac{1}{2}d\bar z_0)+z_2(\del_1-\tfrac{1}{2}d\bar
z_1))
\]
whose bivector component $\beta =(z_0^2-z_1z_2)\del_1\del_2 +
(z_1^2-z_2z_0)\del_2\del_0+(z_2^2-z_0z_1)\del_0\del_1$ is a
quadratic holomorphic Poisson structure on $\C^3$.  This induces a
Poisson structure on $\C P^2$ vanishing on the zero set of the
following cubic polynomial:
\begin{align*}
{\euler}\wedge\beta &= (z_0^3 + z_1^3 + z_2^3 -
3z_0z_1z_2)\del_0\del_1\del_2\\
&=(z_0+z_1+z_2)(z_0+\lambda z_1 +
\lambda^2z_2)(z_0+\lambda^2z_1+\lambda z_2)\del_0\del_1\del_2,
\end{align*}
where $ \euler = \sum z_i \del_i$ is the holomorphic Euler vector field  and $\lambda$ is a third root of unity.  We see that the vanishing
set of this Fermat cubic is the union of three distinct lines in the
plane which intersect at the points $\{[1:1:1],
[1:\lambda:\lambda^2],[1:\lambda^2:\lambda]\}$.

The deformed complex structure can be written explicitly:
\begin{align*}
\gf =e^{\eps}dz_0dz_1dz_2&= (1+\eps)dz_0dz_1dz_2\\
&=(\tfrac{1}{2}(-z_0^2+z_1z_2)dz_0+c.p.)\exp(-\tfrac{1}{2}\tfrac{z_1^2
+ z_0z_2}{-z_2^2 + z_0z_1}dz_1d\bar z_2 +
\tfrac{1}{2}\tfrac{z_0^2+z_1z_2}{-z_2^2+z_0z_1}dz_0d\bar z_2 +
c.p.),
\end{align*}
where ``c.p.'' denotes cyclic permutations of $\{0,1,2\}$. The pure
differential form $\gf$ is of real index zero as long as it has
nonvanishing Mukai pairing with its complex conjugate:
$$
\IP{\gf,\overline\gf} = (\tfrac{R^4}{4} -1) dz_0dz_1dz_2d\bar
z_0d\bar z_1d\bar z_2,
$$
where $R^2 = |z_0|^2 + |z_1|^2 + |z_2|^2$.   The generalized almost
complex structure determined by $\gf$ on the ball of radius $\sqrt
2$ is not integrable, however, since
\[
d\gf = \tfrac{1}{2}(z_0d\bar z_0 + z_1d\bar z_1 + z_2d\bar
z_2)dz_0dz_1dz_2.
\]
Nonetheless, when pulled back to the unit sphere in $\C^3$ this
derivative vanishes, and hence we may proceed as before, quotienting
by the $S^1$-action \eqref{circle action in C3}, as we do next.

\vskip6pt
\noindent
{\bf The reduction.}
We begin with the generalized complex structure $\gf$:
\[
 \gf =(\tfrac{1}{2}(-z_0^2+z_1z_2)dz_0+c.p.)\exp(-\tfrac{1}{2}\tfrac{z_1^2
+ z_0z_2}{-z_2^2 + z_0z_1}dz_1d\bar z_2 +
\tfrac{1}{2}\tfrac{z_0^2+z_1z_2}{-z_2^2+z_0z_1}dz_0d\bar z_2 +
c.p.).
\]
As in Example \ref{example 1},  we calculate the interior product by $\del_\theta$:
\begin{align*}
i_{\del_\theta}\gf=-\tfrac{i(z_0^3+z_1^3+z_2^3-3z_0z_1z_2)}{2}\exp(
&\tfrac{(z_2|z_0|^2 + z_2|z_1|^2 + z_0^2\bar z_1 + 2z_0z_1\bar z_2 +
z_1^2\bar z_0)dz_0dz_1 + (z_1^3-z_2^3)dz_0d\bar z_0
)}{2(z_0^3+z_1^3+z_2^3-3z_0z_1z_2)}\notag\\
&+ \tfrac{(z_0^2z_1-2z_0z_2^2 + z_1^2 z_2)dz_0d\bar z_1 - (z_0^2z_2
- 2 z_0z_1^2 + z_1z_2^2)dz_0d\bar
z_2}{2(z_0^3+z_1^3+z_2^3-3z_0z_1z_2)}+ c.p.)
\end{align*}
Now we pull back to $S^5$ by imposing $1=R^2 = \sum_i z_i\bar z_i$
and obtain a homogeneous differential form after rescaling:
\begin{align*}
\tilde\gf =\exp(&\tfrac{(z_2|z_0|^2 + z_2|z_1|^2 + z_0^2\bar z_1 +
2z_0z_1\bar z_2 + z_1^2\bar z_0)dz_0dz_1 + (z_1^3-z_2^3)dz_0d\bar
z_0}{2R^2(z_0^3+z_1^3+z_2^3-3z_0z_1z_2)}\notag\\
& + \tfrac{(z_0^2z_1-2z_0z_2^2 + z_1^2 z_2)dz_0d\bar z_1 - (z_0^2z_2
- 2 z_0z_1^2 + z_1z_2^2)dz_0d\bar z_2
}{2R^2(z_0^3+z_1^3+z_2^3-3z_0z_1z_2)}+c.p. )
\end{align*}
As in the previous example, we subtract from this the quantity
$\tfrac{dR}{R}\wedge i_{{\euler}+\bar{\euler}}\tilde\gf$, obtaining finally a
manifestly projective representative for the generator for the
canonical bundle:

\begin{align*}
\gf_B=\exp(&\tfrac{((z_1^3-z_2^3 - z_0^3-z_0z_1z_2)|z_1|^2 -
(z_2^3-z_0^3 -z_1^3-z_0z_1z_2)|z_2|^2 +2z_0^2(z_2^2\bar z_1
-z_1^2\bar z_2))dz_0d\bar
z_0}{2R^2(z_0^3+z_1^3+z_2^3-3z_0z_1z_2)}\\
&+\tfrac{(z_1\bar z_0(z_0^3-z_1^3+z_2^3+z_0z_1z_2) -2z_0\bar
z_2(z_1^3+z_2^3-z_0z_1z_2) -2|z_0|^2z_0z_2^2 +
2|z_2|^2z_2z_1^2)dz_0d\bar
z_1}{2R^2(z_0^3+z_1^3+z_2^3-3z_0z_1z_2)}\notag\\
&+ \tfrac{(z_2\bar z_0(-z_0^3-z_1^3 + z_2^3-z_0z_1z_2) +2z_0\bar z_1(z_1^3+z_2^3-z_0z_1z_2) +2|z_0|^2z_0z_1^2-2|z_1|^2z_1z_2^2)dz_0d\bar z_2}{2R^2(z_0^3+z_1^3+z_2^3-3z_0z_1z_2)}\notag\\
&+\tfrac{(z_2(|z_0|^2|z_1|^2 + z_0^2\bar z_1\bar z_2+ \bar z_0^2
z_1z_2 + c.p.))dz_0dz_1}{2R^2(z_0^3+z_1^3+z_2^3-3z_0z_1z_2)}\notag\\
&+c.p.)\notag
\end{align*}
This differential form is closed, but blows up along the three
distinct lines of the type change locus, where one can verify by
rescaling that it defines a complex structure.  This generalized
complex structure, together with the Fubini-Study symplectic
structure, forms a generalized K\"ahler structure on $\C P^2$.

In affine coordinates $(z_1,z_2)$  for $\C P^2$, the type change
locus consists of three lines intersecting at $\{(1,1),
(\lambda,\lambda^2),(\lambda^2,\lambda)\}$. Define $r^2 = z_1\bar
z_1 + z_2\bar z_2$.  We may now write $\gf_B$ in these coordinates:
\begin{align*}
\gf_B =\exp(&\tfrac{((z_2^3-1 - z_1^3-z_1z_2)|z_2|^2 - (1-z_1^3
-z_2^3-z_1z_2) +2z_1^2(\bar z_2 -z_2^2))dz_1d\bar
z_1}{2(1+r^2)(1+z_1^3+z_2^3-3z_1z_2)}\\
 &+\tfrac{((1-z_1^3 -
z_2^3-z_1z_2) - (z_1^3-z_2^3 -1-z_1z_2)|z_1|^2 +2z_2^2(z_1^2 -\bar
z_1))dz_2d\bar z_2}{2(1+r^2)(1+z_1^3+z_2^3-3z_1z_2)}\notag\\
&+\tfrac{(z_2\bar z_1(z_1^3-z_2^3+1+z_1z_2) -2z_1(z_2^3+1-z_1z_2)
-2|z_1|^2z_1 + 2z_2^2)dz_1d\bar
z_2}{2(1+r^2)(1+z_1^3+z_2^3-3z_1z_2)}\notag\\
&+\tfrac{(z_1\bar z_2(-z_2^3-1 + z_1^3-z_1z_2) +2z_2(1+z_1^3-z_1z_2)
+2|z_2|^2z_2-2z_1^2)
dz_2d\bar z_1}{2(1+r^2)(1+z_1^3+z_2^3-3z_1z_2)}\notag\\
&+\tfrac{(|z_1|^2|z_2|^2 + z_1^2\bar z_2+ \bar z_1^2 z_2 +
z_1(|z_2|^2 + z_2^2\bar z_1+ \bar z_2^2z_1) +z_2(|z_1|^2 + \bar
z_1\bar
z_2+z_1z_2))dz_1dz_2}{2(1+r^2)(1+z_1^3+z_2^3-3z_1z_2)})\notag
\end{align*}
This form, together with the Fubini-Study symplectic structure
\begin{align*}
\gf_{A} &= \exp(-\tfrac{1}{2}\tfrac{(1+r^2)(dz_1d\bar z_1 +
dz_2d\bar z_2) -(\bar z_1dz_1+ \bar z_2dz_2)(z_1d\bar z_1 + z_2d\bar
z_2))}{(1+r^2)^2}),
\end{align*}
defines explicitly a generalized K\"ahler structure on $\C P^2$ with
type change along a triangle as described above.
\end{ex}


\vspace{10pt}

\noindent\textsc{Henrique Bursztyn}, Instituto de Matem\'atica Pura
e Aplicada, Estrada Dona Castorina 110, Rio de Janeiro, 22460-320,
Brasil.

\vspace{10pt}

\noindent\textsc{Gil R. Cavalcanti}, Mathematical Institute, 24-29
St. Giles, Oxford, OX1 3LB, UK.

\vspace{10pt}

\noindent\textsc{Marco Gualtieri}, Department of Mathematics,
Massachusetts Institute of Technology, Cambridge, MA 02139, USA.

\end{document}